\documentclass[11pt]{article}
\usepackage[T1]{fontenc}
\usepackage{mathpazo,amsmath,amssymb,amsthm,mathtools}
\usepackage[margin=1in]{geometry}
\usepackage{microtype,booktabs,enumitem}
\usepackage[numbers,sort&compress]{natbib}
\usepackage[hidelinks]{hyperref}
\usepackage[capitalize,nameinlink]{cleveref}
\usepackage{fancyhdr}
\numberwithin{equation}{section}
\newtheorem{theorem}{Theorem}[section]
\newtheorem{proposition}[theorem]{Proposition}
\newtheorem{lemma}[theorem]{Lemma}
\newtheorem{corollary}[theorem]{Corollary}
\theoremstyle{definition}
\newtheorem{assumption}{Assumption}
\crefname{assumption}{Assumption}{Assumptions}
\Crefname{assumption}{Assumption}{Assumptions}
\newtheorem{example}[theorem]{Example}
\theoremstyle{remark}
\newtheorem{remark}[theorem]{Remark}
\newcommand{\R}{\mathbb R}
\newcommand{\E}{\mathbb E}
\newcommand{\eu}{\mathrm{e}}
\newcommand{\U}{U_{\mathrm{ad}}}
\newcommand{\op}{\mathrm{op}}
\newcommand{\F}{\mathrm F}

\newcommand{\tr}{\operatorname{tr}}
\newcommand{\EE}[2]{\mathcal E_h(#1,#2)}
\newcommand{\Mthree}{\mathsf M_3}
\newcommand{\dd}{\,\mathrm d}
\hypersetup{pdftitle={Sampling Error in Time-Discretized Stochastic Optimal Control},pdfauthor={Aurya Javeed and Johannes Milz}}

\title{Sampling Error in Time-Discretized\\Stochastic Optimal Control\thanks{\textbf{Funding:} This work is supported by the Sandia Laboratory Directed Research and Development Program, the National Science Foundation under Award No.~DMS-2410944, and the Ronald J.\ and Carol T.\ Beerman Early Career Professorship.}}
\author{Aurya Javeed\thanks{Optimization and Uncertainty Quantification, Sandia National Laboratories, Albuquerque, NM 87185 (\texttt{asjavee@sandia.gov}).}
\and Johannes Milz\thanks{H.\ Milton Stewart School of Industrial and Systems Engineering, Georgia Institute of Technology, Atlanta, GA 30332 (\texttt{johannes.milz@isye.gatech.edu}).}}
\date{November 23, 2026}
\begin{document}
\maketitle

\begin{abstract}
We study sampling error for Euler discretizations of stochastic optimal control problems driven by Brownian motion.
For ordinary Monte Carlo sampling, the error can diverge under time refinement when empirical-mean fluctuations accumulate as bias.
Under global Lipschitz assumptions and a bounded noise coefficient, we derive finite-sample bounds for the value-function error without assuming smooth value functions or independence across sampling batches.
The bounds isolate the empirical-mean contribution.
For convergence under time refinement, ordinary Monte Carlo generally requires the number of samples per stage to grow at a faster rate than the number of intervals in the time mesh, whereas antithetic and recentered sampling eliminate the empirical-mean term and yield error bounds uniform over meshes with a fixed sample size per stage.
We establish lower bounds showing that the ordinary and antithetic sampling terms are attained.
We also prove sample-size convergence for all three methods on fixed meshes.
For antithetic and recentered sampling, convergence is uniform over meshes in every noise dimension, with sharp sample-size rates in dimensions at least three.
\end{abstract}



\section{Introduction}\label{sec:intro}

Continuous-time stochastic optimal control problems arise in applications such as constrained investment \citep{Zhu2022} and management of biological populations \citep{Tran2022}.
Time discretization is a standard starting point for their numerical solution.
Over a fixed horizon $[0,T]$, time discretization produces a finite sequence of decision stages, and the associated value functions are computed backward through the Bellman recursion.
At each stage, the recursion adds the current cost to the expected continuation value and minimizes the result over the available actions.

The expected continuation value cannot be evaluated analytically in general.
Deterministic quadrature is effective in low noise dimension but becomes expensive as dimension grows, so sampling provides a natural alternative.
In convex settings, for example, stochastic dual dynamic programming and related cutting-plane methods approximate Bellman value functions using sampled expectations \citep{Lan2024,Fuellner2025}, including in time-discretized consumption--investment problems \citep{Pagnoncelli2024}.

We consider problems driven by Brownian motion and a uniform Euler discretization with step size $h=T/N$.
The one-step expectations are then Gaussian. Let $V_k^h(x)$ denote the value function at time $t_k=kh$ and state $x\in\mathbb{R}^n$ when these Gaussian expectations are evaluated exactly, and let $\widehat V_k^h(x)$ denote the value obtained when each expectation is replaced by an average over $M$ sample points.
We refer to these as the \emph{population} and \emph{empirical} value functions, respectively.
At a given stage, the same sample points are used for all states and candidate actions. We assume that the Bellman recursions are otherwise solved exactly, so that the only error is
\(
|\widehat V_k^h(x)-V_k^h(x)|.
\)

Refining the time mesh creates more Bellman updates and therefore more opportunities for sampling error to affect the computation of empirical value functions.
This interaction can be severe.
In \cref{ex:mean}, the Euler discretization itself introduces no value-function error, yet ordinary Monte Carlo sampling with a fixed number of points per stage produces an expected error that diverges as $h\to0$.
The mechanism is simple: the minimization in Bellman recursions converts fluctuations in the empirical means into a systemic error, and that error accumulates over the stages.

\begin{example}[Error accumulation]
\label[example]{ex:mean}
Let $W$ be a one-dimensional standard Brownian motion.
We consider the controlled system
\begin{equation*}
 \dd X_t=u_t\dd W_t,\qquad u_t\in\{-1,1\},
\end{equation*}
with terminal cost $g(x)=x$ and no running cost.
The general formulations of the continuous- and discrete-time control problems are given in
\cref{sec:problem,sec:time-discretization}, respectively.
Under every admissible control, the expected terminal cost
equals the starting state, both in continuous time and for
the Euler discretization with exact expectations. Hence
$V_k^h(x)=x$ for $k=0,\ldots,N$.

At stage $j=1,\ldots,N$, we use the batch
$Z_j^{(1)},\ldots,Z_j^{(M)}$, where all $NM$ samples are
independent with common law $\mathcal N(0,1)$, the standard
normal distribution.
We write
$b_j=M^{-1}\sum_{i=1}^M Z_j^{(i)}$ for the empirical mean
of batch $j$.
The empirical Bellman
recursion is given by
\[
 \widehat V_N^h(x)=x,
 \qquad
 \widehat V_k^h(x)
 =\min_{u\in\{-1,1\}}\frac1M\sum_{i=1}^M
 \widehat V_{k+1}^h\bigl(x+h^{1/2}uZ_{k+1}^{(i)}\bigr),
 \qquad k=0,\ldots,N-1.
\]
At the final step, $\widehat V_{N-1}^h(x)=x-h^{1/2}|b_N|$.
Since the subtracted term is independent of the state,
backward induction gives
\begin{equation*}
 \widehat V_k^h(x)
 =x-h^{1/2}\sum_{j=k+1}^N|b_j|.
\end{equation*}
The empirical mean $b_j$ is distributed according to
$\mathcal N(0,1/M)$, so
\begin{equation*}
 \E[V_0^h(x)-\widehat V_0^h(x)]
 =T\sqrt{2/\pi}\,(hM)^{-1/2},
\end{equation*}
which diverges as $h\to0$ with $M$ fixed.
In contrast, if every batch has zero empirical mean, then
$\widehat V_k^h(x)=x$ for all $k$.
\end{example}

We call a sampling rule \emph{centered} if its empirical
mean vanishes for every realization at every stage.
Two simple constructions have this property. For even $M$,
\emph{antithetic sampling} pairs $M/2$ independent samples
with their negatives, while \emph{recentered sampling}
subtracts the empirical mean from each of the $M$ samples
in a batch. We compare these rules with ordinary Monte
Carlo sampling, which uses $M$ independent samples per stage.

\paragraph{Contributions.}

Let $\Gamma=(\gamma_1,\ldots,\gamma_N)$ be a sequence of sampling rules, let $\widehat V_{k,\Gamma}^h$ denote the corresponding empirical value functions, and let $b_j$ be the mean of $\gamma_j$.
Under global Lipschitz assumptions in the state, a bounded noise coefficient, and Monte Carlo-type bounds on the centered moments of the sampling rules, we prove, for
$h\le1$,
\begin{equation*}
 \E\max_{0\le k\le N}
 \|\widehat V_{k,\Gamma}^h-V_k^h\|_\infty
 \le C\biggl[
 h^{1/2}\sum_{j=1}^N\E|b_j|
 +h^{1/4}+M^{-1/2}
 \biggr].
\end{equation*}
The expectation is over the sampled batches.
The constant is independent of \(h\) and \(M\), the batches may be arbitrarily dependent across stages, and no smoothness of the value functions is assumed.

For ordinary Gaussian sampling, $\E|b_j|$ is of order $M^{-1/2}$.
Since there are $T/h$ stages, the first term is thus of order $(hM)^{-1/2}$.
Antithetic and recentered sampling are centered, so this term vanishes and the bound becomes $C(h^{1/4}+M^{-1/2})$.
The $h^{1/4}$ term is an artifact of smoothing the value functions and is not an intrinsic sampling error.

A separate comparison shows that all three Monte Carlo methods converge as $M\to\infty$ on every fixed mesh.
For antithetic and recentered sampling, the convergence is uniform over time meshes; when the noise dimension $m\ge3$, their error is bounded by
\[
  C\min\{h^{1/4}+M^{-1/2},\,M^{-1/m}\}.
\]

Finally, lower-bound examples show that the $(hM)^{-1/2}$ term for ordinary sampling and the $M^{-1/2}$ term for antithetic sampling are attained at a fixed initial state under time refinement.
A separate example shows that the exponent $1/m$ cannot be improved for bounds uniform over states and time meshes.

\paragraph{Related literature.}

Kushner and Dupuis~\citep{Kushner2001} develop discrete-time approximations of continuous-time stochastic optimal control problems based on local consistency of transition moments.
Krylov~\citep{Krylov1999} restricts the admissible controls of the continuous problem to be piecewise constant on uniform time intervals and obtains an $O(h^{1/6})$ error bound.
The analysis uses ``shaking'', in which small perturbations to the problem appear as additional controls, followed by mollification of the resulting value functions.
The same regularization technique is applied to finite-difference approximations of Bellman equations \citep{Krylov2000}.
Jakobsen et~al.~\citep[Theorem~2.1]{Jakobsen2019} refine this argument and improve the piecewise-constant control error bound to $O(h^{1/4})$.
Their Section~3.2 also considers centered approximating increments with approximate covariance matching and higher-moment conditions tied to the time step.

The work of Picarelli and Reisinger~\citep{Picarelli2020} is closest to ours.
Their scheme combines piecewise-constant controls, Euler time stepping, deterministic Gauss--Hermite quadrature, and spatial interpolation. Their analysis uses Krylov's shaking and mollification technique, which we adapt.
Their quadrature rules satisfy prescribed polynomial moment conditions exactly, whereas our integration rules are random in general and have realization-dependent moment errors.
Our estimates isolate the sampling error relative to the Gaussian Euler recursion.

Using ordinary Gaussian sampling in a Bellman recursion is a sample-average approximation \citep{Kleywegt2002}.
Haskell et~al.~\citep{Haskell2016} analyze empirical value and policy iteration for discounted Markov decision processes.
Their main analysis assumes finite state and action spaces, while their continuous-state inventory example restricts the state to a compact interval.
Fahim et~al.~\citep[Section~4]{Fahim2011} study Monte Carlo sampling under time refinement in a numerical scheme for fully nonlinear parabolic equations.
More recently, Milz and Shapiro~\citep{Milz2025CLT} and Chen et~al.~\citep{Chen2026} study limit theorems for empirical value functions with a fixed number of stages and in discounted infinite-horizon problems, respectively.
Our focus is instead on finite-sample bounds for sampled Euler Bellman recursions.
Our bounds control the expected supremum of the error over states in $\R^n$ and over all stages, with explicit dependence on both $h$ and $M$.

\paragraph{Organization.}
\Cref{sec:main} states the main results, and
\cref{sec:proof} contains their proofs.
\cref{sec:lower} gives the lower-bound and sharpness results.

\section{Setting and main results}\label{sec:main}

\subsection{Problem and assumptions}\label{sec:problem}

Let $T>0$ be the final time. Let $W$ be a standard $m$-dimensional Brownian
motion, and let $\U\subset\R^p$ be nonempty and closed. We
call the elements of $\U$ actions. For $t\in[0,T)$, an admissible
control on $[t,T]$ is a $\U$-valued process $(u_s)$ that is progressively
measurable with respect to the filtration generated by
the increments
$(W_s-W_t)_{s\in[t,T]}$. Starting from $x\in\R^n$ at time
$t$, the state under an admissible control satisfies
\begin{equation}\label{eq:continuous-sde}
\dd X_s=f(s,X_s,u_s)\dd s+\sigma(s,X_s,u_s)\dd W_s,
\qquad X_t=x.
\end{equation}
Given running and terminal costs $\ell$ and $g$, respectively,
the continuous-time value function is the optimal expected
remaining cost:
\begin{equation*}
 V(t,x)=\inf_u\E\left[\int_t^T\ell(s,X_s,u_s)\dd s+g(X_T)\right],
\end{equation*}
where the infimum is over admissible controls on $[t,T]$.
At the final time we set $V(T,x)=g(x)$.
Background on continuous-time stochastic control is given in \citep[Chapters~III--IV]{Fleming2006}.

We write $|\cdot|$ for the Euclidean norm on vectors,
$\|\cdot\|_{\op}$ for the induced matrix operator norm, and
$\|\cdot\|_{\F}$ for the Frobenius norm on matrices and
tensors. Let $I_m$ denote the identity matrix in
$\R^{m\times m}$. If $z\in\R^m$, then $z^{\otimes3}$ has
entries $z_a z_b z_c$ and
$\|z^{\otimes3}\|_{\F}=|z|^3$. For functions
$\phi:\R^n\to\R$, we set
$\|\phi\|_\infty=\sup_{x\in\R^n}|\phi(x)|$.

\begin{assumption}[Problem data]\label[assumption]{ass:data}
The functions $f:[0,T]\times\R^n\times\U\to\R^n$, $\sigma:[0,T]\times\R^n\times\U\to\R^{n\times m}$, $\ell:[0,T]\times\R^n\times\U\to\R$, and $g:\R^n\to\R$ are Borel measurable. The functions $f,\sigma,\ell$ are continuous in $u$
for each $(t,x)$. There exist $L_1,L_2<\infty$ such that, uniformly in $t,x,y,u$,
\begin{equation*}
\begin{aligned}
 |f(t,x,u)-f(t,y,u)|+\|\sigma(t,x,u)-\sigma(t,y,u)\|_{\op}&\le L_1|x-y|,\\
 |f(t,0,u)|+\|\sigma(t,x,u)\|_{\op}&\le L_1,\\
 |g(x)-g(y)|+|\ell(t,x,u)-\ell(t,y,u)|&\le L_2|x-y|,\\
 |g(0)|+|\ell(t,0,u)|&\le L_2.
\end{aligned}
\end{equation*}
\end{assumption}
The Lipschitz and growth conditions in \cref{ass:data} are
uniform in time and action. They allow linear growth of the
drift and costs, while requiring the noise coefficient to be
bounded, and $\U$ need not be compact. These conditions ensure
a unique strong solution to the state equation for each
admissible control.
Time continuity is not assumed.

Throughout the paper, $C$ denotes an upper-bound constant depending
only on $T,n,m,L_1,L_2$, unless further dependencies are
indicated, and may change between displays.
We use $C_m$ for
constants depending only on $m$ and $c>0$ for lower-bound constants.
All constants are
independent of $h$ and $M$.

\subsection{Time discretization}\label{sec:time-discretization}

The time mesh is given by $h=T/N$ and $t_k=kh$, $k=0,\ldots,N$.
We assume $h\le1$. For $z\in\R^m$, we define
\begin{equation*}
 F_k^h(x,u,z)
 =x+hf(t_k,x,u)+h^{1/2}\sigma(t_k,x,u)z,
 \qquad k=0,\ldots,N-1.
\end{equation*}
The Euler discretization of \eqref{eq:continuous-sde}
is the successive application of this transition,
substituting independent random vectors with law
$P=\mathcal N(0,I_m)$ for the argument $z$.

To compare Gaussian expectations with sampling rules, we allow
a sequence of noise laws $\Lambda=(\lambda_1,\ldots,\lambda_N)$
in place of $P$, each with a finite first moment.
Let $Z_j\sim\lambda_j$, $j=1,\ldots,N$, and let
$(\mathcal H_j)_{j=0}^N$ be a filtration such that
$Z_j$ is $\mathcal H_j$-measurable and $Z_{j+1}$ has
conditional law $\lambda_{j+1}$ given $\mathcal H_j$, written
\begin{equation}\label{eq:conditional-law}
 \mathcal L(Z_{j+1}\mid\mathcal H_j)=\lambda_{j+1}.
\end{equation}
We allow $(\mathcal H_j)$ to contain auxiliary information
beyond the past noise increments.

For a starting stage
$k=0,\ldots,N-1$, a control $u=(u_j)_{j=k}^{N-1}$ is admissible if
each $u_j$ is $\mathcal H_j$-measurable and takes values
in $\U$.
Given such a control, let $(S_j)_{j=k}^N$ denote the
discrete state process. Then
\begin{equation*}
 S_{j+1}=F_j^h(S_j,u_j,Z_{j+1}),
 \qquad S_k=x,
 \qquad j=k,\ldots,N-1.
\end{equation*}
The corresponding value functions are
\begin{equation}\label{eq:discrete-cost}
 v_N^\Lambda=g,\qquad
 v_k^\Lambda(x)
 \coloneqq\inf_u\E\biggl[
 h\sum_{j=k}^{N-1}\ell(t_j,S_j,u_j)+g(S_N)
 \biggr],
 \qquad k=0,\ldots,N-1,
\end{equation}
where the infimum is over admissible controls.
We suppress the dependence on $h$ in this notation.

For Borel functions
$\phi:\R^n\to\R$ of at most linear growth, we define the Bellman operator
$T_k^\Lambda$ by
\begin{equation}\label{eq:bellman}
 [T_k^\Lambda\phi](x)
 =\inf_{u\in\U}\Bigl\{
 h\ell(t_k,x,u)
 +\int\phi(F_k^h(x,u,z))\lambda_{k+1}(\mathrm d z)
 \Bigr\}.
\end{equation}
The Bellman recursion is then given by
\begin{equation}\label{eq:recursion}
 v_N^\Lambda=g,\qquad
 v_k^\Lambda=T_k^\Lambda v_{k+1}^\Lambda,
 \qquad k=0,\ldots,N-1.
\end{equation}
Under \cref{ass:data}, \cref{lem:wellposed}\ref{item:wellposed-policies}
shows that this recursion gives the value functions in
\eqref{eq:discrete-cost} and that the infimum there is unchanged when
restricted to controls of the form $u_j=\pi_j(S_j)$ for
Borel maps $\pi_j:\R^n\to\U$. We call such controls
\emph{Markov policies}. An analogous Bellman recursion and Markov policy
representation with Gaussian Euler increments and a compact action set
appear in \citep[Theorem~3.4 and Remark~3.7]{Bonnans2019}. A proof sketch
under the present assumptions follows \cref{lem:wellposed}.
The computation of $v_k^\Lambda$ from
$v_{k+1}^\Lambda$ is called a \emph{Bellman update}.
A single law in the superscript means that law
at every stage; thus the population value functions
from \cref{sec:intro} satisfy $V_k^h=v_k^P$.

Under \cref{ass:data}, the value functions $v_k^\Lambda$,
$k=0,\ldots,N$, are finite, of at most linear growth,
and globally Lipschitz by
\cref{lem:wellposed}\ref{item:wellposed-regularity}.
The proof uses the elementary inequality
\begin{equation}\label{eq:infimum-comparison}
 \Bigl|\inf_{a\in A}\phi(a)-\inf_{a\in A}\psi(a)\Bigr|
 \le\sup_{a\in A}|\phi(a)-\psi(a)|.
\end{equation}
The functions $\phi,\psi$ are real-valued on a common
nonempty set $A$, and both infima are assumed finite.
The growth and Lipschitz constants may depend on $h$
and the sequence $\Lambda$.

For two sequences of noise laws, $\Lambda=(\lambda_j)$ and
$\Lambda'=(\lambda_j')$, we use the shorthand
\begin{equation*}
 \EE{\Lambda}{\Lambda'}
 =\max_{0\le k\le N}
 \|v_k^\Lambda-v_k^{\Lambda'}\|_\infty.
\end{equation*}

\subsection{Sampling}\label{sec:rules}

A \emph{sampling rule} is a random, finitely supported
probability measure on $\R^m$. We use one for each
stage:
\begin{equation*}
 \gamma_j=\sum_i\omega_{ji}\delta_{\eta_{ji}},
 \qquad \omega_{ji}\ge0,\qquad \sum_i\omega_{ji}=1,
 \qquad j=1,\ldots,N.
\end{equation*}
For each $j$, the sum contains a deterministic finite number
of terms. The weights $\omega_{ji}$ and points $\eta_{ji}$
are measurable random variables, and points may repeat.
The notation $\delta_z$ denotes unit mass at $z$.
We call the collection at one stage a batch.
We write $\Gamma=(\gamma_1,\ldots,\gamma_N)$.

Substituting $\Gamma$ for $\Lambda$ in
\eqref{eq:recursion} defines
the empirical value functions
$\widehat V_{k,\Gamma}^h=v_k^\Gamma$,
$k=0,\ldots,N$.
With $\Lambda=\Gamma$, the integrals in
\eqref{eq:bellman} are finite weighted sums:
\[
 \int\phi(F_k^h(x,u,z))\,\gamma_{k+1}(\mathrm d z)
 =
 \sum_i\omega_{k+1,i}
 \phi(F_k^h(x,u,\eta_{k+1,i})),
 \qquad k=0,\ldots,N-1.
\]
In \eqref{eq:discrete-cost}, all batches are held fixed.
Each $Z_j\sim\gamma_j$ corresponds
to selecting an entry from the $j$th batch according
to its weight.
The expectation averages over trajectories formed
by making these selections independently across stages.
The batches themselves may be arbitrarily dependent
across stages, e.g., the same batch may be
reused at every stage.

We write $b_j$ and $\mu_j$, respectively, for the mean and
centered version of $\gamma_j$:
\begin{equation}\label{eq:center}
 b_j=\int z\,\gamma_j(\mathrm d z),\qquad
 \mu_j=\sum_i\omega_{ji}\delta_{\eta_{ji}-b_j},
\qquad j=1,\ldots,N.
\end{equation}
A sampling rule $\gamma_j$ is \emph{centered} if
$b_j=0$ almost surely.
The centered moments of the sampling
rule are
\begin{equation}\label{eq:centered-moments}
 Q_j=\int zz^\top\mu_j(\mathrm d z),\qquad
 H_j=\int z^{\otimes3}\mu_j(\mathrm d z),\qquad
 R_j=\int|z|^4\mu_j(\mathrm d z).
\end{equation}
We impose the following assumption on these centered moments.
In contrast to the expectation in \eqref{eq:discrete-cost},
$\E$ in the assumption averages over the random batches.
\begin{assumption}[Moment errors]
\label[assumption]{ass:sampling}
For an accuracy parameter $M\ge1$, we assume that
\begin{equation}\label{eq:moment-hyp}
 \E\|Q_j-I_m\|_{\F}^2\le K/M,\qquad
 \E\|H_j\|_{\F}\le KM^{-1/2},\qquad
 \E R_j\le K,
 \qquad j=1,\ldots,N,
\end{equation}
where $K$ is an accuracy parameter independent of $h,M,j$.
\end{assumption}
For the Gaussian sampling rules below, $M$ counts the points
in each batch with multiplicity, as in \cref{ex:mean} for ordinary sampling.
For a general sampling rule, $M$ need not equal the number of terms in a batch.
Our first main result is a finite-sample bound for the expected error between the empirical and population value functions.

\begin{theorem}[Finite-sample error]
\label[theorem]{thm:main}
Under \textup{\cref{ass:data,ass:sampling}},
\begin{equation}\label{eq:main-bound}
 \E\EE{\Gamma}{P}
 \le C_K\biggl[
 h^{1/2}\sum_{j=1}^N\E|b_j|
 +h^{1/4}+M^{-1/2}
 \biggr].
\end{equation}
Here $C_K$ depends only on $T,n,m,L_1,L_2,K$.
\end{theorem}

The proof of \cref{thm:main} is given in
\cref{sec:assembly}.
For centered sampling rules, the empirical-mean term in
\eqref{eq:main-bound} vanishes. The moment conditions
\eqref{eq:moment-hyp} are still required, however.

\subsection{Ordinary, antithetic, and recentered Monte Carlo}
\label{sec:montec}

We apply \cref{thm:main} to three Gaussian sampling rules.
For each stage $j=1,\ldots,N$, let $Z_j^{(1)}, \ldots, Z_j^{(M)}$ be independent vectors
with law $P$. Batches at different stages
may be arbitrarily dependent. The ordinary, antithetic,
and recentered Monte Carlo sampling rules are, respectively,
\begin{equation*}
 \gamma_j^{\mathcal O}
 =\frac1M\sum_{i=1}^M\delta_{Z_j^{(i)}},
\qquad
 \gamma_j^{\mathcal A}
 =\frac1{M}\sum_{i=1}^{M/2}
   \bigl(\delta_{Z_j^{(i)}}+\delta_{-Z_j^{(i)}}\bigr),
\qquad
 \gamma_j^{\mathcal R}
 =\frac1M\sum_{i=1}^M\delta_{Z_j^{(i)}-\bar Z_j},
   \qquad \bar Z_j=\frac1M\sum_{i=1}^M Z_j^{(i)}.
\end{equation*}
Each rule uses $M$ points, counted with multiplicity.
Antithetic sampling requires $M$ to be even and pairs
each of $M/2$ independent Gaussian draws with its negative.
Both the antithetic and recentered rules are centered.
Symmetry also gives $H_j=0$
for the antithetic rule, whereas the recentered rule
need not have a vanishing third centered moment.

For $\mathcal S\in\{\mathcal O,\mathcal A,\mathcal R\}$,
we write $\Gamma^{\mathcal S}
=(\gamma_1^{\mathcal S},\ldots,\gamma_N^{\mathcal S})$.
The Gaussian moment calculations in \cref{app:moments} show
that all three rules satisfy \cref{ass:sampling} with $K$
depending only on $m$. Hence \cref{thm:main} implies the
following bounds.

\begin{corollary}[Gaussian sampling]\label[corollary]{cor:samplers}
Under \textup{\cref{ass:data}}, ordinary Monte Carlo sampling for $M\ge1$ satisfies
\begin{equation*}
 \E\EE{\Gamma^{\mathcal O}}{P}\le C(h^{1/4}+(hM)^{-1/2}).
\end{equation*}
Antithetic sampling for even $M\ge2$ and recentered sampling for $M\ge1$ satisfy
\begin{equation}\label{eq:centered-bound}
 \E\EE{\Gamma^{\mathcal S}}{P}\le C(h^{1/4}+M^{-1/2}),
 \qquad \mathcal S\in\{\mathcal A,\mathcal R\}.
\end{equation}
\end{corollary}

The term $(hM)^{-1/2}$ is attained in \cref{ex:mean},
and the term $M^{-1/2}$ is attained by antithetic sampling in
\cref{sec:cosine-example}. On a fixed time mesh, the term $h^{1/4}$ in
\cref{cor:samplers} does not vanish as $M\to\infty$.
Our second main result shows that the sampling error
nevertheless vanishes as $M\to\infty$.

\begin{theorem}[Consistency]
\label[theorem]{thm:transport}
Under \textup{\cref{ass:data}}, for each
$\mathcal S\in\{\mathcal O,\mathcal A,\mathcal R\}$,
we have $\E\EE{\Gamma^{\mathcal S}}{P}\to0$
as $M\to\infty$ on each fixed mesh.
For antithetic and recentered sampling, this convergence
is uniform over $h=T/N\le1$.
Furthermore, if $m\ge3$, antithetic and recentered sampling satisfy
\begin{equation}\label{eq:centered-combined}
 \E\EE{\Gamma^{\mathcal S}}{P}
 \le C\min\{h^{1/4}+M^{-1/2},\,M^{-1/m}\},
 \qquad \mathcal S\in\{\mathcal A,\mathcal R\}.
\end{equation}
The constant is independent of $h$ and $M$, and all
conclusions allow arbitrary dependence between stage
batches. The sample-size restrictions are those in
\textup{\cref{cor:samplers}}.
\end{theorem}

The proof is given in \cref{sec:transport}.
The $1/m$ exponent is sharp for bounds uniform over
meshes in the global supremum norm as shown in
\cref{prop:unweighted-obstruction}.

\subsection{Continuous-time consequences}\label{sec:continuous}

We now show that our main results combine with
time-discretization
error estimates to characterize total error and sufficient
conditions for convergence to the continuous-time
value function $V$.

Let $V^h$ and $\widehat V_\Gamma^h$ denote the
piecewise constant extensions in time of the population
and empirical value functions, respectively:
\[
 V^h(t,x)=V_k^h(x),\qquad
 \widehat V_\Gamma^h(t,x)=\widehat V_{k,\Gamma}^h(x),
 \qquad t\in[t_k,t_{k+1}),
 \qquad k=0,\ldots,N-1.
\]
Both extensions equal $g(x)$ at $t=T$.
To bound the error between $V$ and $\widehat V_\Gamma^h$,
we use $V^h$ as an intermediate approximation.
We define the population approximation error by
\begin{equation*}
 \mathcal D(h)=\sup_{t\in[0,T]}
 \left\|
 \frac{V^h(t,\cdot)-V(t,\cdot)}
      {(1+|\cdot|^2)^{1/2}}
 \right\|_\infty.
\end{equation*}
The denominator of this expression accounts for
time-discretization errors that grow linearly in $|x|$.
Such errors can occur even under \cref{ass:data}.
For example, we consider a problem with $n=m=1$, no
control, $f(x)=x$, $\sigma=\ell=0$, and $g(x)=x$.
Then $V(0,x)=\eu^Tx$ and $V^h(0,x)=(1+h)^Nx$.
For every $h>0$, the absolute difference between
these functions grows without bound as $|x|\to\infty$.

The triangle inequality, $(1+|x|^2)^{-1/2}\le1$,
and \cref{thm:main} give the following total error estimate.

\begin{corollary}[Total error]\label[corollary]{cor:total}
Under the hypotheses of \textup{\cref{thm:main}},
\[
 \E\sup_{t\in[0,T]}
 \left\|
 \frac{\widehat V_\Gamma^h(t,\cdot)-V(t,\cdot)}
      {\sqrt{1+|\cdot|^2}}
 \right\|_\infty
 \le \mathcal D(h)+C_K\biggl[
 h^{1/2}\sum_{j=1}^N\E|b_j|
 +h^{1/4}+M^{-1/2}
 \biggr].
\]
If $\mathcal D(h)=O(h^{1/4})$, the respective total bounds are
\begin{equation*}
 C(h^{1/4}+(hM)^{-1/2})\quad\text{and}\quad C_K(h^{1/4}+M^{-1/2})
\end{equation*}
for ordinary Monte Carlo sampling and any centered sampling
rule satisfying \eqref{eq:moment-hyp}. If only $\mathcal D(h)\to0$
is known, sufficient conditions for total consistency are
$h\to0$ and $hM\to\infty$ in the ordinary case, and
$h\to0$ and $M\to\infty$ in the centered case.
\end{corollary}

The assumptions in \cref{ass:data} alone do not ensure
$\mathcal D(h)\to0$, since they impose no time regularity.
\Cref{prop:population} gives sufficient conditions for
an unweighted population error bound of order $h^{1/4}$.
Under those conditions, the corresponding total error
bounds hold without the weight.

\section{Proofs and sampling bounds}\label{sec:proof}

We use two arguments to compare value functions
associated with different sequences of noise laws.
The trajectory comparison in \cref{sec:strong}
couples controlled Euler trajectories using the same
controls. It bounds the value function error in terms
of noise means and Wasserstein distances between
centered laws.
The second argument compares Bellman updates for laws
with matching means and covariances. For smooth
functions, Taylor expansion terms through second order
cancel in the comparison. Shaking and mollification
provide smooth approximations of value functions satisfying a one-sided
Bellman inequality (\cref{sec:shaking}).
This inequality allows the update bounds to be
propagated backward (\cref{sec:weak}).
We combine these comparisons for fixed sequences of
noise laws in \cref{sec:assembly}, then average over
the random sampling batches to prove \cref{thm:main}.
The trajectory comparison also yields the consistency
and sampling-rate results in \cref{sec:transport}.

\subsection{Trajectory comparison}
\label{sec:strong}

For the trajectory comparison, we use the Wasserstein 2-distance, defined for probability laws $\lambda,\mu$ with finite
second moments by
\begin{equation*}
 W_2^2(\lambda,\mu)
 =\inf\{\E|Z-Y|^2:
          \mathcal L(Z)=\lambda,\ \mathcal L(Y)=\mu\}.
\end{equation*}
The infimum is over all joint distributions with the specified
marginals, called couplings.

\begin{lemma}[Trajectory comparison]\label[lemma]{lem:strong}
Under \textup{\cref{ass:data}}, for $j=1,\ldots,N$, let $\lambda_j$
be a centered noise law satisfying
$\int |z|^2\lambda_j(\mathrm d z)\le \kappa$, and let $\gamma_j$
be a noise law with finite second moment, mean $b_j$, and
centered law $\mu_j$. We set $\Lambda=(\lambda_1,\ldots,\lambda_N)$ and
$\Gamma=(\gamma_1,\ldots,\gamma_N)$. Then
\begin{equation}\label{eq:strong}
 \EE{\Gamma}{\Lambda}
 \le C_\kappa\biggl[h^{1/2}\sum_{j=1}^N|b_j|
       +\Bigl(h\sum_{j=1}^N W_2^2(\mu_j,\lambda_j)\Bigr)^{1/2}\biggr].
\end{equation}
The constant $C_\kappa$ depends only on $T,n,m,L_1,L_2,\kappa$.
\end{lemma}
\begin{proof}
At $k=N$, $v_N^\Lambda=v_N^\Gamma=g$.
For $k=0,\ldots,N-1$ and $x\in\R^n$, we choose couplings
$(\xi_j,\zeta_j)$ of $\lambda_j,\mu_j$,
$j=k+1,\ldots,N$, such that
\[
  d_j^2=\E|\xi_j-\zeta_j|^2
\]
is arbitrarily close to $W_2^2(\lambda_j,\mu_j)$.
We take these pairs independently across stages.
All expectations in this proof are taken with respect to
these coupled increments. The laws are fixed.

We take $\mathcal H_k$ to be the trivial $\sigma$-algebra and set
\[
 \mathcal H_j=\sigma(\xi_i,\zeta_i:k<i\le j),
 \qquad j=k+1,\ldots,N.
\]
The filtration $(\mathcal H_j)_{j=k}^N$ records the joint
history of the coupled increments.
We fix any control $u=(u_j)_{j=k}^{N-1}$ admissible with
respect to $(\mathcal H_j)_{j=k}^N$.
Starting from $S_k=\widetilde S_k=x$, we construct the
coupled Euler trajectories
\begin{align*}
 S_{j+1}
 &=S_j+hf(t_j,S_j,u_j)
       +h^{1/2}\sigma(t_j,S_j,u_j)\xi_{j+1},\\
 \widetilde S_{j+1}
 &=\widetilde S_j+hf(t_j,\widetilde S_j,u_j)
       +h^{1/2}\sigma(t_j,\widetilde S_j,u_j)
                    (\zeta_{j+1}+b_{j+1}),
 \qquad j=k,\ldots,N-1.
\end{align*}
The growth bounds in \cref{ass:data} and the second-moment
assumptions on the laws imply by induction that
$S_j$ and $\widetilde S_j$ have finite second moments at every stage.

We first estimate the state difference, separating the contribution
of the means $b_j$. We write
\[
 \delta S_j=S_j-\widetilde S_j,
\qquad
 \delta f_j=f(t_j,S_j,u_j)-f(t_j,\widetilde S_j,u_j),
\qquad
 \delta\sigma_j
 =\sigma(t_j,S_j,u_j)-\sigma(t_j,\widetilde S_j,u_j).
\]
Subtracting the state equations gives
\begin{equation}\label{eq:strong-recursion}
 \delta S_{j+1}
 =\delta S_j+h\delta f_j+h^{1/2}B_{j+1}
       -h^{1/2}\sigma(t_j,\widetilde S_j,u_j)b_{j+1},
\end{equation}
where
\[
 B_{j+1}
 =\delta\sigma_j\xi_{j+1}
  +\sigma(t_j,\widetilde S_j,u_j)
       (\xi_{j+1}-\zeta_{j+1}).
\]

We estimate the right-hand side of \eqref{eq:strong-recursion}
in $L^2$ with respect to the joint law of the coupled increments,
treating the contribution of $b_{j+1}$ separately. We write
$\|X\|_{L^2}=(\E|X|^2)^{1/2}$.
The term containing $b_{j+1}$ satisfies
\begin{equation}\label{eq:mean-contribution}
 \bigl\|h^{1/2}\sigma(t_j,\widetilde S_j,u_j)b_{j+1}
       \bigr\|_{L^2}
 \le L_1h^{1/2}|b_{j+1}|.
\end{equation}
To estimate $\delta S_j+h\delta f_j+h^{1/2}B_{j+1}$
in $L^2$, we evaluate the conditional moments of $B_{j+1}$.
The pair $(\xi_{j+1},\zeta_{j+1})$ is independent of
$\mathcal H_j$ and both components are centered.
Since the coefficients in $B_{j+1}$ are
$\mathcal H_j$-measurable,
\[
 \E[B_{j+1}\mid\mathcal H_j]=0,
\]
while the elementary inequality $|v+w|^2\le2|v|^2+2|w|^2$ and the bounds
on $\sigma$ in \cref{ass:data} also give
\[
 \E[|B_{j+1}|^2\mid\mathcal H_j]
 \le 2L_1^2\kappa|\delta S_j|^2+2L_1^2d_{j+1}^2.
\]

The zero conditional mean makes $B_{j+1}$ orthogonal in
$L^2$ to the $\mathcal H_j$-measurable vector
$\delta S_j+h\delta f_j$.
The Pythagorean theorem therefore gives
\begin{equation}\label{eq:centered-contribution}
 \|\delta S_j+h\delta f_j+h^{1/2}B_{j+1}\|_{L^2}^2
 =\|\delta S_j+h\delta f_j\|_{L^2}^2
   +h\|B_{j+1}\|_{L^2}^2
 \le (1+C_\kappa h)\|\delta S_j\|_{L^2}^2
   +Chd_{j+1}^2.
\end{equation}
The inequality follows from taking expectations in the conditional
second-moment bound for $B_{j+1}$ and using
$|\delta f_j|\le L_1|\delta S_j|$ as well as $h\le1$.
Applying the triangle inequality in $L^2$ to
\eqref{eq:strong-recursion} and using
\eqref{eq:mean-contribution} and \eqref{eq:centered-contribution}
yields
\begin{equation}\label{eq:l2-iteration}
 \|\delta S_{j+1}\|_{L^2}
 \le
 \Bigl((1+C_\kappa h)\|\delta S_j\|_{L^2}^2
             +Chd_{j+1}^2\Bigr)^{1/2}
       +L_1h^{1/2}|b_{j+1}|.
\end{equation}

We iterate \eqref{eq:l2-iteration}. When all $b_i$ vanish,
the steps that follow reduce to the usual discrete Gronwall argument
for $\|\delta S_j\|_{L^2}^2$.
For nonzero $b_i$, we use the triangle inequality in $\mathbb R^2$
to keep their accumulated contribution outside the square root.
Applying \eqref{eq:l2-iteration} at stage $j-1$ gives a bound
for $\|\delta S_j\|_{L^2}$. Substituting this bound into
the square-root term of \eqref{eq:l2-iteration} at stage $j$
gives the following estimate for some $\alpha\ge0$:
\[
 \Bigl((1+C_\kappa h)
        [\alpha+L_1h^{1/2}|b_j|]^2+Chd_{j+1}^2\Bigr)^{1/2}
 \le
 (1+C_\kappa h)^{1/2}L_1h^{1/2}|b_j|
 +\Bigl((1+C_\kappa h)\alpha^2+Chd_{j+1}^2\Bigr)^{1/2}.
\]
Repeating this step and using $\delta S_k=0$ gives by induction
\[
 \begin{aligned}
 \|\delta S_j\|_{L^2}
 \le L_1h^{1/2}\sum_{i=k+1}^j
       (1+C_\kappa h)^{(j-i)/2}|b_i|
 +
 \Bigl(Ch\sum_{i=k+1}^j
       (1+C_\kappa h)^{j-i}d_i^2\Bigr)^{1/2},
 \quad j=k,\ldots,N,
 \end{aligned}
\]
and since $(1+C_\kappa h)^N\le \eu^{C_\kappa T}$,
\begin{equation}\label{eq:state-strong}
 \max_{k\le j\le N}\|\delta S_j\|_{L^2}
 \le C_\kappa\biggl[
       h^{1/2}\sum_{i=k+1}^N|b_i|
       +\Bigl(h\sum_{i=k+1}^N d_i^2\Bigr)^{1/2}
       \biggr].
\end{equation}
This bound is uniform over the initial state $x$ and
the common adapted control $u$.

We now estimate the value function
error. By \eqref{eq:discrete-cost}, $v_k^\Lambda(x)$ and
$v_k^\Gamma(x)$ are
infima over the same class of controls admissible with respect
to $(\mathcal H_j)_{j=k}^N$.
Applying \eqref{eq:infimum-comparison}, Lipschitz continuity
of the costs, and the Cauchy--Schwarz inequality gives
\[
 |v_k^\Lambda(x)-v_k^\Gamma(x)|
 \le L_2\sup_u\biggl(
       h\sum_{j=k}^{N-1}\|\delta S_j\|_{L^2}
       +\|\delta S_N\|_{L^2}
       \biggr)\\
 \le C_\kappa\biggl[
       h^{1/2}\sum_{i=k+1}^N|b_i|
       +\Bigl(h\sum_{i=k+1}^N d_i^2\Bigr)^{1/2}
       \biggr].
\]
The second inequality follows from \eqref{eq:state-strong}
and $(N-k)h\le T$.
We let each $d_i^2$ tend to $W_2^2(\lambda_i,\mu_i)$ and
extend the sums to all stages.
Taking the supremum over initial states and the maximum
over starting stages proves \eqref{eq:strong}.
\end{proof}

\subsection{Shaking and mollification}\label{sec:shaking}

We construct smooth approximations to the value functions that
satisfy a one-sided form of the original Bellman recursion.
For state-dependent coefficients, mollifying the value functions
directly does not in general preserve this inequality.
We adapt the argument in
\citep[Propositions~4.1--4.2]{Picarelli2020} to centered noise
laws with uniformly bounded second moments.

\begin{lemma}[Uniform regularity]
\label[lemma]{lem:regularity}
Under \textup{\cref{ass:data}}, let $\Lambda=(\lambda_1,\ldots,\lambda_N)$
be a sequence of centered noise laws satisfying
\(
 m_2(\lambda_j):=\int |z|^2\lambda_j(\mathrm d z)\le\kappa<\infty,
\)
\(
 j=1,\ldots,N.
\)
Then, for $k=0,\ldots,N$ and $x,y\in\R^n$,
\begin{equation}\label{eq:regularity}
 |v_k^\Lambda(x)|\le C_\kappa(1+|x|),\qquad
 |v_k^\Lambda(x)-v_k^\Lambda(y)|\le C_\kappa|x-y|.
\end{equation}
The constant depends only on $T,n,m,L_1,L_2,\kappa$.
The same estimates hold when the actions are enlarged to
$(u,e)$ with $u\in\U$ and $|e|\le\varepsilon$, and $f,\sigma,\ell$
are evaluated at $x+e$, while the terminal cost remains $g$.
The constants are uniform for $0<\varepsilon\le1$.
\end{lemma}

\begin{proof}
For $k=N$, the bounds follow directly from $v_N^\Lambda=g$.
For $k=0,\ldots,N-1$ and $x,y\in\R^n$, we construct
trajectories with $S_k=x$ and $\widetilde S_k=y$,
using the same noise increments with laws
$\lambda_{k+1},\ldots,\lambda_N$ and the same
admissible controls.
For $j=k,\ldots,N-1$,
\[
 \E|S_{j+1}|^2
 \le (1+C_\kappa h)\E|S_j|^2+C_\kappa h,
\qquad
 \E|S_{j+1}-\widetilde S_{j+1}|^2
 \le (1+C_\kappa h)\E|S_j-\widetilde S_j|^2.
\]
The first inequality follows from $\lambda_{j+1}$ being
centered, $m_2(\lambda_{j+1})\le\kappa$, and
\cref{ass:data}.
The second follows from the calculation leading to
\eqref{eq:l2-iteration} with $b_i=d_i=0$.
Discrete Gronwall arguments imply
\[
 \max_{k\le j\le N}\E|S_j|^2
 \le C_\kappa(1+|x|^2),
\qquad
 \max_{k\le j\le N}\E|S_j-\widetilde S_j|^2
 \le C_\kappa|x-y|^2.
\]
The growth and Lipschitz bounds for the costs, together
with the Cauchy--Schwarz inequality and
\eqref{eq:infimum-comparison}, now give \eqref{eq:regularity} in a manner similar to the last display of the proof of \Cref{lem:strong}.

The shifted coefficients and running cost satisfy the
bounds in \cref{ass:data} with constants independent of
$0<\varepsilon\le1$. The same argument therefore applies
to the enlarged controls $(u,e)$.
\end{proof}

\begin{lemma}[Smooth approximation]
\label[lemma]{lem:shaking}
Under the hypotheses of \textup{\cref{lem:regularity}}, for each
$0<\varepsilon\le1$ there exist
$w_0,\ldots,w_N\in C^\infty(\R^n)$ such that
\begin{align}
  \max_{0\le k\le N}\|w_k-v_k^\Lambda\|_\infty
  &\le C_\kappa\varepsilon,
       \label{eq:smooth-close}\\
  \max_{0\le k\le N}\|D^rw_k\|_\infty
  &\le C_{\kappa,r}\varepsilon^{1-r},
       \qquad r=1,2,\ldots,
       \label{eq:smooth-derivs}\\
  w_k&\le T_k^\Lambda w_{k+1},
       \qquad k=0,\ldots,N-1.
       \label{eq:smooth-bellman}
\end{align}
Here $\|\cdot\|_{\mathrm{op}}$ denotes the multilinear
operator norm, and
$\|D^rw_k\|_\infty
=\sup_{x\in\R^n}\|D^rw_k(x)\|_{\mathrm{op}}$.
The constants depend only on $T,n,m,L_1,L_2,\kappa$
and the fixed mollifier, with additional dependence on $r$
for $C_{\kappa,r}$. No higher-moment bounds are required.
\end{lemma}

\begin{proof}
We perturb the control problem by allowing small spatial
shifts in the coefficients and running cost. In the resulting
shaken problem, the actions are pairs $(u,e)$ with
$u\in\U$ and $|e|\le\varepsilon$, and at state $x\in\R^n$ the functions
$f,\sigma,\ell$ are evaluated at $x+e$.
Let $v_0^{\Lambda,\varepsilon},\ldots,
v_N^{\Lambda,\varepsilon}$ be the value functions of this
problem, defined as in \eqref{eq:discrete-cost} with terminal
cost $g$.
The transitions and Bellman operators are, for $k = 0,\ldots,N-1$,
\[
\begin{aligned}
  F_k^{h,e}(x,u,z)
  &=x+hf(t_k,x+e,u)+h^{1/2}\sigma(t_k,x+e,u)z,\\
  [T_k^{\Lambda,\varepsilon}\phi](x)
  &=\inf_{\substack{u\in\U\\|e|\le\varepsilon}}
    \biggl\{
      h\ell(t_k,x+e,u)
      +\int\phi(F_k^{h,e}(x,u,z))
             \lambda_{k+1}(\mathrm d z)
    \biggr\}.
\end{aligned}
\]
The shifted data remain continuous in $(u,e)$, so
\cref{lem:wellposed} ensures that the shaken value functions
satisfy the Bellman recursion with these operators.

We next compare the shaken value functions
with the original value functions $v_0^\Lambda, \ldots, v_N^\Lambda$.
For $k=0,\ldots,N-1$ and any admissible control
$(u_j,e_j)_{j=k}^{N-1}$, we couple the original trajectory
$(S_k,\ldots,S_N)$ and the shaken trajectory
$(\widetilde S_k,\ldots,\widetilde S_N)$, starting from
$S_k=\widetilde S_k=x$ and taking, at each stage
$j=k,\ldots,N-1$, the same action $u_j$ and the same
noise increment with law $\lambda_{j+1}$.
We write $\delta S_j=S_j-\widetilde S_j$.
The differences in $f$ and in $\sigma$ are bounded by
$L_1(|\delta S_j|+\varepsilon)$. The conditional
second-moment calculation leading to
\eqref{eq:l2-iteration} therefore gives
\[
  \E|\delta S_{j+1}|^2
  \le (1+C_\kappa h)\E|\delta S_j|^2
       +C_\kappa h\varepsilon^2.
\]
A discrete Gronwall argument with $\delta S_k=0$ implies
\[
 \max_{k\le j\le N}\E|\delta S_j|^2
 \le C_\kappa\varepsilon^2.
\]

Both families of value functions may be written as infima
over the enlarged class of admissible controls
$(u_j,e_j)_{j=k}^{N-1}$.
Applying \eqref{eq:infimum-comparison}, the Lipschitz bounds
in \cref{ass:data}, the Cauchy--Schwarz inequality, and
$(N-k)h\le T$ results in
\begin{equation}\label{eq:shaken-close}
  |v_k^\Lambda(x)-v_k^{\Lambda,\varepsilon}(x)| \le L_2\sup_{u,e}\biggl[
       h\sum_{j=k}^{N-1}
         \bigl(\E|\delta S_j|+\varepsilon\bigr)
       +\E|\delta S_N|
       \biggr]
  \le C_\kappa\varepsilon.
\end{equation}
Since $v^\Lambda_N = v^{\Lambda,\varepsilon}_N = g$, this estimate also holds for $k=N$.

We now mollify the shaken value functions.
We choose a nonnegative function $\varphi\in C^\infty(\R^n)$,
supported in the unit ball with integral one, and set
\[
 \varphi_\varepsilon(y)=\varepsilon^{-n}\varphi(y/\varepsilon),
 \qquad
 w_k(x)=\int v_k^{\Lambda,\varepsilon}(x-y)
                 \varphi_\varepsilon(y)\dd y,
  \qquad k=0,\ldots,N.
\]
The growth bound in \cref{lem:regularity} and the smoothness
and compact support of $\varphi_\varepsilon$ ensure that
$w_k$ is well defined and smooth.
The same lemma gives
\[
  |w_k(x)-v_k^{\Lambda,\varepsilon}(x)|
  \le \int
       \bigl|v_k^{\Lambda,\varepsilon}(x-y)
             -v_k^{\Lambda,\varepsilon}(x)\bigr|
       \varphi_\varepsilon(y)\dd y\\
  \le C_\kappa\int |y|\varphi_\varepsilon(y)\dd y
   \le C_\kappa\varepsilon.
\]
This together with \eqref{eq:shaken-close} proves
\eqref{eq:smooth-close}.

For $r\ge1$, integration by parts and compact support implies
$\int D^r\varphi_\varepsilon(y)\dd y=0$.
Thus,
\[
 \|D^rw_k(x)\|_{\mathrm{op}}
  =\Bigl\|\int
       \bigl[v_k^{\Lambda,\varepsilon}(x-y)
             -v_k^{\Lambda,\varepsilon}(x)\bigr]
       D^r\varphi_\varepsilon(y)\dd y\Bigr\|_{\mathrm{op}}
  \le C_\kappa\int
       |y|\,\|D^r\varphi_\varepsilon(y)\|_{\mathrm{op}}
\dd y
  \le C_{\kappa,r}\varepsilon^{1-r},
\]
where the last inequality follows from substituting $y=\varepsilon z$.
The result \eqref{eq:smooth-derivs} follows.

Finally, for $k=0,\ldots,N-1$ and any action $u$, we choose
$e=y$ at state $x-y$ in the shaken recursion.
This $e$ is admissible on the support of $\varphi_\varepsilon$,
and
\(
  F_k^{h,y}(x-y,u,z)=F_k^h(x,u,z)-y.
\)
From the Bellman recursion for the shaken value functions then,
\[
  v_k^{\Lambda,\varepsilon}(x-y)
  \le h\ell(t_k,x,u)
     +\int
       v_{k+1}^{\Lambda,\varepsilon}(F_k^h(x,u,z)-y)
       \lambda_{k+1}(\mathrm d z).
\]
We multiply by $\varphi_\varepsilon(y)$ and integrate over $y$.
The growth bound in \eqref{eq:regularity}, the compact support
of $\varphi_\varepsilon$, and the finite first moment of
$\lambda_{k+1}$ justify exchanging the order of integration.
Consequently,
\[
  w_k(x)
  \le h\ell(t_k,x,u)
     +\int w_{k+1}(F_k^h(x,u,z))
            \lambda_{k+1}(\mathrm d z).
\]
Taking the infimum over $u$ proves \eqref{eq:smooth-bellman}.
\end{proof}

\subsection{Bellman update comparison}
\label{sec:weak}

We compare Bellman updates for centered noise laws with
matching covariances, using the smooth approximations from
\cref{lem:shaking}. Taylor expanding bounds the difference
in terms of the third-moment discrepancy and the fourth
moments. We then propagate this bound backward using
\eqref{eq:smooth-bellman} to estimate the value function error.
For a probability law $\lambda$ with finite fourth moment,
we write
\[
  \Mthree(\lambda)=\int z^{\otimes3}\lambda(\mathrm d z),
  \qquad
  m_4(\lambda)=\int |z|^4\lambda(\mathrm d z).
\]

\begin{proposition}[Matching covariances]
\label[proposition]{prop:weak}
Under \textup{\cref{ass:data}}, let
$\Lambda=(\lambda_1,\ldots,\lambda_N)$ and
$\Lambda'=(\lambda_1',\ldots,\lambda_N')$
be sequences of centered noise laws with finite fourth moments.
Suppose that, for finite $\kappa\ge0$,
\begin{equation}\label{eq:matching-cov}
  \int zz^\top\lambda_j(\mathrm d z)
  =\int zz^\top\lambda_j'(\mathrm d z)
  =\Sigma_j,
  \qquad
  0\preceq\Sigma_j\preceq\kappa I_m,
  \quad j=1,\ldots,N.
\end{equation}
We set $\Delta_j=\Mthree(\lambda_j)-\Mthree(\lambda_j')$.
Then, for every $0<\varepsilon\le1$,
\begin{equation}\label{eq:weak}
  \EE{\Lambda}{\Lambda'}
  \le C_\kappa\biggl[
       \varepsilon
       +h^{3/2}\varepsilon^{-2}
          \sum_{j=1}^N\|\Delta_j\|_{\F}
       +h^2\varepsilon^{-3}
          \sum_{j=1}^N
            \bigl(m_4(\lambda_j)+m_4(\lambda_j')\bigr)
       \biggr].
\end{equation}
The constant $C_\kappa$ depends only on
$T,n,m,L_1,L_2,\kappa$ and the fixed mollifier.
\end{proposition}

\begin{proof}
By \eqref{eq:matching-cov},
$m_2(\lambda_j)=\tr\Sigma_j\le m\kappa$
for $j=1,\ldots,N$.
For $0<\varepsilon\le1$, we apply
\cref{lem:shaking} to obtain smooth approximations
$w_0,\ldots,w_N$ for the $\Lambda$-problem.
For $k=0,\ldots,N-1$, $x\in\R^n$, and $u\in\U$, we write
$z_0=x+hf(t_k,x,u)$ and $A=\sigma(t_k,x,u)$.
Taylor expanding in the noise argument $z$,
\[
\begin{aligned}
  w_{k+1}(F_k^h(x,u,z))
  &=w_{k+1}(z_0)
    +h^{1/2}Dw_{k+1}(z_0)[Az]\\
  &\quad+\frac h2D^2w_{k+1}(z_0)[Az,Az]
  +\frac{h^{3/2}}6
       D^3w_{k+1}(z_0)[Az,Az,Az]
       +\mathcal R(z).
\end{aligned}
\]
By \eqref{eq:smooth-derivs} and
$\|A\|_{\mathrm{op}}\le L_1$,
\[
  |\mathcal R(z)|
  \le \frac{h^2}{24}\|D^4w_{k+1}\|_\infty|Az|^4
  \le C_\kappa h^2\varepsilon^{-3}|z|^4.
\]

We integrate the expansion of $w_{k+1}$ against the laws
$\lambda_{k+1}$ and $\lambda_{k+1}'$ to bound the difference between the
associated Bellman updates.
Since both laws are centered and have the same covariance,
the terms through second order cancel in the difference.
Writing $A_a$ for column $a$ of $A$, the difference of
the cubic terms is
\[
  \frac{h^{3/2}}6
  \sum_{a,b,c=1}^m
    D^3w_{k+1}(z_0)[A_a,A_b,A_c]\,
    (\Delta_{k+1})_{abc}.
\]
Each derivative coefficient is bounded in absolute value
by $C_\kappa\varepsilon^{-2}$, while, from the Cauchy-Schwarz inequality,
\[
  \sum_{a,b,c=1}^m|(\Delta_{k+1})_{abc}|
  \le m^{3/2}\|\Delta_{k+1}\|_{\F}.
\]
Consequently, the definition of the Bellman operators and
\eqref{eq:infimum-comparison} imply
\[
\begin{aligned}
  \bigl|[T_k^\Lambda w_{k+1}](x)
        -[T_k^{\Lambda'}w_{k+1}](x)\bigr|
  &\le \sup_{u\in\U}
       \Bigl|\int w_{k+1}(F_k^h(x,u,z))\,
         (\lambda_{k+1}-\lambda_{k+1}')(\mathrm d z)\Bigr|\\
  &\le C_\kappa\biggl[
       h^{3/2}\varepsilon^{-2}\|\Delta_{k+1}\|_{\F}
       +h^2\varepsilon^{-3}
          \bigl(m_4(\lambda_{k+1})
                +m_4(\lambda_{k+1}')\bigr)
       \biggr]
       =:r_{k+1}.
\end{aligned}
\]

Combining this result with \eqref{eq:smooth-bellman} gives
\[
  w_k
  \le T_k^\Lambda w_{k+1}
  \le T_k^{\Lambda'}w_{k+1}+r_{k+1}.
\]
We propagate this bound backward.
Each $T_k^{\Lambda'}$ is monotone and satisfies $T_k^{\Lambda'}(\phi+c)=T_k^{\Lambda'}\phi+c$ for constants $c$.
Hence a pointwise bound
$w_{k+1}\le v_{k+1}^{\Lambda'}+c$ implies
\[
  w_k
  \le T_k^{\Lambda'}(v_{k+1}^{\Lambda'}+c)+r_{k+1}
  =v_k^{\Lambda'}+c+r_{k+1}.
\]
Since $w_N\le g+C_\kappa\varepsilon =v_N^{\Lambda'}+C_\kappa\varepsilon$,
backward induction and \eqref{eq:smooth-close} therefore establish
\[
  v_k^\Lambda
  \le w_k+C_\kappa\varepsilon
  \le v_k^{\Lambda'}+2C_\kappa\varepsilon
       +\sum_{j=k+1}^{N}r_j,
  \qquad k=0,\ldots,N-1.
\]
To obtain the inequality with $v_k^\Lambda$ and
$v_k^{\Lambda'}$ interchanged, we construct smooth
approximations for the $\Lambda'$-problem and repeat
the argument.
We extend the sums to all stages and take the supremum
over $x$ and the maximum over $k$.
Substituting the definition of $r_j$ then proves
\eqref{eq:weak}.
\end{proof}

\begin{remark}[Higher moments]
\label[remark]{rem:remainder}
In the setting of \cref{prop:weak}, let $r\ge3$
be an integer. Suppose that the laws match all
polynomial moments through degree $r$ at each
stage and have uniformly bounded $(r+1)$st
absolute moments.
Taylor expanding through degree $r$ and choosing
$\varepsilon=h^{1/2-1/(r+1)}$ give
$\EE{\Lambda}{\Lambda'}\le C_rh^{1/2-1/(r+1)}$.
Thus matching third moments gives the rate $h^{1/4}$,
while higher orders give exponents approaching $1/2$.
The constants depend on $r$ and the moment bound.
\end{remark}

\subsection{Combining comparisons}
\label{sec:assembly}

We now combine the trajectory and Bellman update
comparisons to prove \cref{thm:main}.
For a fixed sequence of sampling rules
$\Gamma=(\gamma_1,\ldots,\gamma_N)$, we construct
two intermediate sequences: $\widetilde\Gamma$,
obtained by rescaling the centered rules so that their
covariance matrices are bounded above by $I_m$, and
$G$, consisting of centered Gaussian laws with the
same covariances as $\widetilde\Gamma$ at each stage.
The triangle inequality gives
\begin{equation}\label{eq:three-comparisons}
  \EE{\Gamma}{P}
  \le \EE{\Gamma}{\widetilde\Gamma}
     +\EE{\widetilde\Gamma}{G}
     +\EE{G}{P}.
\end{equation}
We bound the first and third terms using
\cref{lem:strong}, and the middle term using
\cref{prop:weak} with $\kappa=1$.

Applied to the centered version of each $\gamma_j$,
the following lemma constructs $\widetilde\Gamma$ and $G$
and provides estimates needed
to bound the right-hand side of
\eqref{eq:three-comparisons}.
The proof of the lemma consists of routine calculations and is deferred to \cref{app:normalization}.

\begin{lemma}[Covariance normalization]
\label[lemma]{lem:normalization}
Let $\mu$ be a centered law on $\R^m$ with covariance
$Q$ and finite fourth moment. We set
$a=\max\{1,\|Q\|_{\mathrm{op}}^{1/2}\}$ and
$\widetilde Q=Q/a^2$.
For $Z\sim\mu$, we define
\[
  \widetilde\mu=\mathcal L(Z/a),
  \qquad
  G_Q=\mathcal N(0,\widetilde Q).
\]
Both of these laws are centered with covariance $\widetilde Q$,
where $0\preceq\widetilde Q\preceq I_m$, and
\begin{equation}\label{eq:normalization-w2}
  W_2(\mu,\widetilde\mu)+W_2(G_Q,P)
  \le C_m\|Q-I_m\|_{\F}.
\end{equation}
The constant $C_m$ depends only on $m$.
Moreover,
\[
\Mthree(\widetilde\mu)=a^{-3}\Mthree(\mu),
\qquad
m_4(\widetilde\mu)\le m_4(\mu),
\qquad
m_4(G_Q)\le m(m+2).
\]
\end{lemma}

We now bound the right-hand side of
\eqref{eq:three-comparisons} in terms of the moments
of $\Gamma$.

\begin{proposition}[Deterministic error bound]
\label[proposition]{prop:deterministic}
Under \textup{\cref{ass:data}}, let
$\Gamma=(\gamma_1,\ldots,\gamma_N)$ be a sequence
of sampling rules, and let $P=\mathcal N(0,I_m)$.
The quantities $b_j,Q_j,H_j,R_j$ are defined in
\eqref{eq:center}--\eqref{eq:centered-moments}.
There is a constant $C$, depending only on
$T,n,m,L_1,L_2$ and a fixed mollifier, such that,
for every $0<\varepsilon\le1$,
\[
  \EE{\Gamma}{P}\le C\biggl[
    h^{1/2}\sum_{j=1}^N|b_j|
    +\Bigl(h\sum_{j=1}^N
       \|Q_j-I_m\|_{\F}^2\Bigr)^{1/2}
    +\varepsilon
    +h^{3/2}\varepsilon^{-2}
       \sum_{j=1}^N\|H_j\|_{\F}
    +h^2\varepsilon^{-3}
       \sum_{j=1}^N(1+R_j)
    \biggr].
\]
In particular, choosing $\varepsilon=h^{1/4}$ implies
\begin{equation}\label{eq:deterministic}
  \EE{\Gamma}{P}\le C\biggl[
    h^{1/2}\sum_{j=1}^N|b_j|
    +\Bigl(h\sum_{j=1}^N
       \|Q_j-I_m\|_{\F}^2\Bigr)^{1/2}
    +h\sum_{j=1}^N\|H_j\|_{\F}
    +h^{1/4}\Bigl(1+h\sum_{j=1}^N R_j\Bigr)
    \biggr].
\end{equation}
\end{proposition}
\begin{proof}
For $j=1,\ldots,N$, let $\mu_j$ be the centered
version of $\gamma_j$ defined in \eqref{eq:center}.
We apply \cref{lem:normalization} to $\mu_j$ to obtain
$a_j,\widetilde\mu_j,\widetilde Q_j$ and $G_{Q_j}$.
We set $\widetilde\Gamma=(\widetilde\mu_1,\ldots,\widetilde\mu_N)$
and $G=(G_{Q_1},\ldots,G_{Q_N})$.

The laws $\widetilde\mu_j$ and $P$ are centered with
second moments at most $m$.
Applying \cref{lem:strong} with these reference laws
to the first and third terms of
\eqref{eq:three-comparisons}, respectively, and using
\eqref{eq:normalization-w2}, results in
\[
  \EE{\Gamma}{\widetilde\Gamma}+\EE{G}{P}
  \le C\biggl[
    h^{1/2}\sum_{j=1}^N|b_j|
    +\Bigl(h\sum_{j=1}^N
       \|Q_j-I_m\|_{\F}^2\Bigr)^{1/2}
    \biggr].
\]

For the middle term of \eqref{eq:three-comparisons},
$\widetilde\mu_j$ and $G_{Q_j}$ are centered with
common covariance $0\preceq\widetilde Q_j\preceq I_m$.
The Gaussian third moments vanish, so
$\Delta_j=a_j^{-3}H_j$ and
$\|\Delta_j\|_{\F}\le\|H_j\|_{\F}$.
Moreover, \cref{lem:normalization} gives
$m_4(\widetilde\mu_j)+m_4(G_{Q_j})\le R_j+m(m+2)$.
Thus \cref{prop:weak} with $\kappa=1$ yields
\[
  \EE{\widetilde\Gamma}{G}
  \le C\biggl[
    \varepsilon
    +h^{3/2}\varepsilon^{-2}
       \sum_{j=1}^N\|H_j\|_{\F}
    +h^2\varepsilon^{-3}
       \sum_{j=1}^N(1+R_j)
    \biggr].
\]

Combining these two bounds in
\eqref{eq:three-comparisons} proves the first bound
in the proposition.
Since $h\le1$, we may take $\varepsilon=h^{1/4}$.
Using $Nh=T$ gives \eqref{eq:deterministic}.
\end{proof}

Since the constant in \eqref{eq:deterministic} is
independent of the realized rules, we can average this
estimate over the random points and any random weights
used to construct them.
In the proof below, $\E$ denotes this expectation.

\begin{proof}[Proof of \textup{\cref{thm:main}}]
Jensen's inequality, \cref{ass:sampling}, and $Nh=T$
give
\begin{align*}
  \E\Bigl(h\sum_{j=1}^N
     \|Q_j-I_m\|_{\F}^2\Bigr)^{1/2}
  &\le
    \Bigl(h\sum_{j=1}^N
       \E\|Q_j-I_m\|_{\F}^2\Bigr)^{1/2}
   \le (TK)^{1/2}M^{-1/2},\\
  h\sum_{j=1}^N\E\|H_j\|_{\F}
  &\le TKM^{-1/2},\\
  h^{1/4}\Bigl(1+h\sum_{j=1}^N\E R_j\Bigr)
  &\le h^{1/4}(1+TK).
\end{align*}
Taking expectations in \eqref{eq:deterministic}
and substituting these bounds proves
\eqref{eq:main-bound}.
Dependence between the sampling batches can be arbitrary.
\end{proof}

\subsection{Consistency}
\label{sec:transport}

We prove \cref{thm:transport} by applying
\cref{lem:strong} to the Gaussian sampling rules.
Centering and symmetrization reduce the required
Wasserstein bounds to those for ordinary sampling.

\begin{proof}[Proof of \textup{\cref{thm:transport}}]
Let $\gamma$ be a probability law on $\R^m$ with
finite second moment, mean $b$, and centered law $\mu$.
Subtracting $b$ from the first component of a coupling
of $\gamma$ and $P$ gives a coupling of $\mu$ and $P$.
Since $P$ is centered, this reduces the quadratic
transport cost by $|b|^2$.
The correspondence is reversible, so taking infima gives
\begin{equation}\label{eq:transport-centering}
  W_2^2(\gamma,P)=|b|^2+W_2^2(\mu,P).
\end{equation}
Let $\gamma^-$ be the law of $-Z$ when $Z\sim\gamma$.
Convexity of $W_2^2(\cdot,P)$
\citep[Theorem~4.8]{Villani2009} and symmetry of $P$
give
\begin{equation}\label{eq:transport-symmetrization}
  W_2^2\bigl(\tfrac12(\gamma+\gamma^-),P\bigr)
  \le \tfrac12 W_2^2(\gamma,P)
     +\tfrac12 W_2^2(\gamma^-,P)
  =W_2^2(\gamma,P).
\end{equation}

At a single stage, $\gamma^{\mathcal O}$,
$\gamma^{\mathcal A}$, and $\gamma^{\mathcal R}$
denote the ordinary, antithetic, and recentered
Gaussian sampling rules.
In the remainder of the proof, $\E$ denotes
expectation over the sampling batches.
The ordinary mean satisfies
$\E|b_j|\le C_mM^{-1/2}$, while the antithetic
and recentered means vanish.
Applying \cref{lem:strong} with reference $P$
and averaging over the batches gives
\begin{align}
  \E\EE{\Gamma^{\mathcal O}}{P}
  &\le C\biggl[
    (hM)^{-1/2}
    +\left(\E W_2^2(\gamma^{\mathcal O},P)\right)^{1/2}
    \biggr],
  \label{eq:transport-ordinary}\\
  \E\EE{\Gamma^{\mathcal S}}{P}
  &\le C\left(
    \E W_2^2(\gamma^{\mathcal S},P)
    \right)^{1/2},
  \qquad \mathcal S\in\{\mathcal A,\mathcal R\}.
  \label{eq:transport-centered}
\end{align}
The constant is independent of $h,M$, and these bounds
allow arbitrary dependence between sampling batches.

Since $P$ has moments of every order, the moment bound
in \citep[Theorem~1]{Fournier2015} gives
\[
  \E W_2^2(\gamma^{\mathcal O},P)\to0
\]
as $M\to\infty$.
The recentered rule is obtained by centering an
ordinary Gaussian rule with $M$ points, while the
antithetic rule symmetrizes one with $M/2$ points.
Thus \eqref{eq:transport-centering} and
\eqref{eq:transport-symmetrization} imply the same
convergence for $\gamma^{\mathcal R}$ and
$\gamma^{\mathcal A}$.
Consequently, \eqref{eq:transport-ordinary} gives
consistency of ordinary sampling at each fixed mesh,
while \eqref{eq:transport-centered} gives consistency
of the centered rules uniformly over $h=T/N\le1$,
since its right-hand side does not depend on $h$.

When $m\ge3$, the Gaussian matching bound in \citep[Theorem~1.1]{LedouxZhu2021} gives
\[
  \E W_2^2(\gamma^{\mathcal O},P)\le C_mM^{-2/m}.
\]
The centering and symmetrization inequalities give
the same rate for both centered rules.
Substituting into \eqref{eq:transport-centered} yields
\begin{equation}\label{eq:transport-sharp}
  \sup_{h=T/N\le1}\E\EE{\Gamma^{\mathcal S}}{P}
  \le CM^{-1/m},
  \qquad \mathcal S\in\{\mathcal A,\mathcal R\}.
\end{equation}
Combining \eqref{eq:transport-sharp} with
\eqref{eq:centered-bound} proves
\eqref{eq:centered-combined}.
\end{proof}

More generally, centered rules satisfying
$\sup_{h,j}\E W_2^2(\gamma_j,P)\to0$ as $M\to\infty$
are consistent uniformly over time meshes.
The moment conditions in \cref{ass:sampling} alone
do not ensure consistency at a fixed mesh.

\begin{remark}[Moment conditions and consistency]
\label[remark]{rem:moment-only}
For $M=6q$, $q=1,2,\ldots$, the equal-weight rule
with $4q$ zeros and $q$ copies of each of $\pm\sqrt3$
is the three-point Gauss--Hermite rule
\[
  \gamma_M=\frac16\delta_{-\sqrt3}
          +\frac23\delta_0
          +\frac16\delta_{\sqrt3}.
\]
This rule satisfies \cref{ass:sampling} with $K=3$.
We consider the problem with $T=h=1$, $n=m=1$, no control,
$f=\ell=0$, $\sigma=1$, and $g(x)=|x|$.
\Cref{ass:data} holds, yet
\[
  v_0^P(0)-v_0^{\gamma_M}(0)
  =\sqrt{2/\pi}-1/\sqrt3
  \not\to0.
\]
\end{remark}

\section{Lower bounds and sharpness}\label{sec:lower}

We establish two complementary sharpness results.
The first concerns sampling error at a fixed initial state
under time refinement. We obtain lower bounds of order
$(hM)^{-1/2}$ and $M^{-1/2}$ for the expected errors of
ordinary Monte Carlo and antithetic sampling, respectively.
\Cref{sec:cosine-example} develops the antithetic lower bound, and
\cref{sec:lower-both} combines it with the mean-error construction
from \cref{ex:mean} in a single control problem.
Every admissible control is optimal when the Gaussian expectations are evaluated exactly.
The lower bounds therefore reflect error in value-function
estimation rather than policy suboptimality.
The second result concerns the dimension-dependent lower
bound $M^{-1/m}$.
\Cref{sec:global-sharpness} establishes this bound using a problem with
one time step and no control, and a single bounded
Lipschitz terminal cost fixed independently of both
the sampling rule and $M$.
All examples satisfy \cref{ass:data}, and their discrete-time population
value functions agree with the continuous-time value functions on the time mesh.
Throughout this section, $T=1$,
$h=1/N$, and $t_k=kh$.

\subsection{Sampling error with zero empirical means}
\label{sec:cosine-example}

We consider the scalar control problem driven by a two-dimensional standard Brownian motion $W=(W_1,W_2)$,
\begin{equation}\label{eq:scalar-cosine}
 \dd X_t=u_t\dd W_{1,t}+(1-u_t)\dd W_{2,t},
 \qquad u_t\in\{0,1\},
\end{equation}
with $\ell=0$ and $g(x)=1-\cos x$.
The population Bellman recursion is
\[
 V_N^h(x)=1-\cos x,\quad
 V_k^h(x)=\min_{u\in\{0,1\}}\int_{\R^2}
 V_{k+1}^h\bigl(x+h^{1/2}\,[uz_1+(1-u)z_2]\bigr)P(\mathrm d z),
\quad
k=0,\ldots,N-1.
\]
For $Z\sim\mathcal N(0,1)$, the Gaussian characteristic function gives
\begin{equation}\label{eq:cosine-gaussian}
 \E\cos(x+h^{1/2} Z)
 =\cos x\,\E\cos(h^{1/2} Z)-\sin x\,\E\sin(h^{1/2} Z)
 =\rho_h\cos x,
\qquad
\rho_h:=\eu^{-h/2}.
\end{equation}
It follows from backward induction through the population Bellman recursion that
\begin{equation}\label{eq:scalar-value}
 V_k^h(x)=1-\rho_h^{N-k}\cos x.
\end{equation}

Let $M=2\nu$. For each stage $j=1,\ldots,N$, we take the batch
$Z_j^{(1)},\ldots,Z_j^{(\nu)}\sim P$, with all $N\nu$ samples
independent, and form $\Gamma^{\mathcal A} = (\gamma_1^{\mathcal A},\ldots \gamma_N^{\mathcal A})$ as in \cref{sec:montec}.
We write $Z_{j,r}^{(i)}$ for the $r$th coordinate of $Z_j^{(i)}$.
We define the empirical cosine averages by
\begin{equation*}
 J_{j,r}:=\int_{\R^2}\cos(h^{1/2}\,z_r)\gamma_j^{\mathcal A}(\mathrm d z)
       =\frac1\nu\sum_{i=1}^{\nu}\cos(h^{1/2} Z_{j,r}^{(i)}),
 \qquad j=1,\ldots,N,\qquad r=1,2.
\end{equation*}
At stage $j-1$, actions $u=1$ and $u=0$ correspond to
$J_{j,1}$ and $J_{j,2}$, respectively. The variables $J_{j,r}$
are independent and identically distributed over $j$ and $r$.
A calculation analogous to \eqref{eq:cosine-gaussian} shows
\begin{equation}\label{eq:antithetic-cosine-update}
 \int_{\R^2}\cos(x+h^{1/2}\,z_r)\gamma_j^{\mathcal A}(\mathrm d z)
 =J_{j,r}\cos x,\qquad j=1,\ldots,N,\qquad r=1,2.
\end{equation}
As a result, $\widehat V_{N-1,\Gamma^{\mathcal A}}^h(0)=1-\max\{J_{N,1},J_{N,2}\}$.
Using $\max(a,b)=(a+b+|a-b|)/2$, the bias of this quantity is
\begin{equation*}
 \E[V_{N-1}^h(0)-\widehat V_{N-1,\Gamma^{\mathcal A}}^h(0)]
 =\E\max\{J_{N,1},J_{N,2}\}-\rho_h
 =\tfrac12\E|J_{N,1}-J_{N,2}|
 =:\delta_{h,\nu}.
\end{equation*}

\begin{lemma}[Bounds on $\delta_{h,\nu}$]\label[lemma]{lem:cosine}
For $0<h\le1$ and every integer $\nu\ge1$,
\begin{equation*}
 ch\nu^{-1/2}\le\delta_{h,\nu}\le\tfrac12h\nu^{-1/2},
 \qquad c=\frac{(1-\eu^{-1})^3}{2\sqrt{27}}>0.
\end{equation*}
\end{lemma}
\begin{proof}
We set $\Xi=J_{N,1}-J_{N,2}$. The calculations in \cref{app:cosine} give
\begin{equation}\label{eq:D-moments}
 \E \Xi^2=\frac{(1-\eu^{-h})^2}{\nu},\qquad
 \E \Xi^4\le\frac{27h^4}{\nu^2}.
\end{equation}
The upper bound on $\delta_{h,\nu}$ follows from the Cauchy--Schwarz inequality and $1-\eu^{-h}\le h$, that is,
\[
 \delta_{h,\nu}\le\tfrac12(\E \Xi^2)^{1/2}
 =\tfrac12(1-\eu^{-h})\nu^{-1/2}\le\tfrac12h\nu^{-1/2}.
\]
For the lower bound, H\"older's inequality with exponents $3/2$ and $3$ implies
\[
 \E \Xi^2=\E\bigl[|\Xi|^{2/3}|\Xi|^{4/3}\bigr]
 \le(\E|\Xi|)^{2/3}(\E \Xi^4)^{1/3}.
\]
The concavity of $s\mapsto1-\eu^{-s}$ ensures $1-\eu^{-h}\ge h(1-\eu^{-1})$,
so rearranging the preceding display and using \eqref{eq:D-moments} yields
\[
 2\delta_{h,\nu} = \E|\Xi|
 \ge\frac{(\E \Xi^2)^{3/2}}{(\E \Xi^4)^{1/2}}
 \ge\frac{(1-\eu^{-1})^3h}{\sqrt{27}}\,\nu^{-1/2}.
 \qedhere
\]
\end{proof}

We now show how this bias accumulates over stages.
For each realization of the batches, backward induction
starting from
$\widehat V_{N,\Gamma^{\mathcal A}}^h(x)=1-\cos x$
and using \eqref{eq:antithetic-cosine-update} gives
\[
\begin{aligned}
 \widehat V_{k,\Gamma^{\mathcal A}}^h(x)
 &=\min_{r\in\{1,2\}}
 \int_{\R^2}
 \widehat V_{k+1,\Gamma^{\mathcal A}}^h(x+h^{1/2}z_r)
 \,\gamma_{k+1}^{\mathcal A}(\mathrm d z)\\
 &\le
 \min_{r\in\{1,2\}}
 \left\{
 1-J_{k+1,r}\cos x
 \prod_{j=k+2}^N\max\{J_{j,1},J_{j,2}\}
 \right\}\\
 &\le
 1-\cos x\prod_{j=k+1}^N\max\{J_{j,1},J_{j,2}\},
 \qquad k=N-1,\ldots,0.
\end{aligned}
\]
The first inequality uses the induction hypothesis and
\eqref{eq:antithetic-cosine-update}. The second evaluates
the minimum at an index attaining
$\max\{J_{k+1,1},J_{k+1,2}\}$.
Averaging over the independent batches, we have
\[
 \E\prod_{j=1}^N\max\{J_{j,1},J_{j,2}\}
 =(\rho_h+\delta_{h,\nu})^N.
\]
Hence
\begin{equation}\label{eq:scalar-lower}
 \E[V_0^h(0)-\widehat V_{0,\Gamma^{\mathcal A}}^h(0)]
 \ge(\rho_h+\delta_{h,\nu})^N-\rho_h^N
 \ge N\rho_h^{N-1}\delta_{h,\nu}
 \ge c_*M^{-1/2},
\end{equation}
where the last inequality uses \cref{lem:cosine} with
$Nh=1$, $\rho_h^{N-1}\ge \eu^{-1/2}$, and
\begin{equation}\label{eq:lower-constant}
 c_*:=\sqrt{\frac{2}{\eu}}\,c
 =\frac{(1-\eu^{-1})^3}{\sqrt{54\eu}}>0.
\end{equation}

\subsection{One problem attaining both sampling bounds}\label{sec:lower-both}

We combine the two scalar control problems into a single
two-dimensional problem. With $n=m=2$ and
$\U=\{-1,1\}\times\{0,1\}$, we consider
\begin{equation}\label{eq:lower-model}
 \begin{aligned}
 \dd X_{1,t}&=2^{-1/2}u_{1,t}(\dd W_{1,t}+\dd W_{2,t}),\\
 \dd X_{2,t}&=u_{2,t}\dd W_{1,t}+(1-u_{2,t})\dd W_{2,t}.
 \end{aligned}
\end{equation}
There is no running cost, and the terminal cost is
\[
 g(x_1,x_2)=x_1+1-\cos x_2.
\]
The first state equation reproduces \cref{ex:mean},
driven by the standard Brownian motion
$(W_1+W_2)/\sqrt2$, while the second reproduces
\eqref{eq:scalar-cosine}.

The covariance matrix associated with the noise coefficient is
\begin{equation*}
 \sigma(u)\sigma(u)^\top=
 \begin{pmatrix}
 1&u_1/\sqrt2\\
 u_1/\sqrt2&1
 \end{pmatrix}
 \succeq(1-2^{-1/2})I_2.
\end{equation*}
The problem is therefore uniformly elliptic: its covariance
matrix is bounded below by a positive multiple of the
identity, uniformly over all actions.

\begin{theorem}[Lower bounds for ordinary and antithetic sampling]\label[theorem]{thm:lower}
For the control problem \eqref{eq:lower-model}, the continuous-time value function and the Euler-discretized population value functions satisfy
\begin{equation}\label{eq:lower-values}
 V(t,x)=x_1+1-\eu^{-(1-t)/2}\cos x_2,\qquad V_k^h(x)=V(t_k,x).
\end{equation}
With independent stage batches, ordinary sampling for every integer $M\ge1$ and antithetic sampling for every even integer $M\ge2$ satisfy, respectively,
\begin{align}
 \E[V_0^h(0)-\widehat V_{0,\Gamma^{\mathcal O}}^h(0)]
 &\ge\sqrt{2/\pi}\,(hM)^{-1/2},
 \label{eq:lower-ordinary}\\
 \E[V_0^h(0)-\widehat V_{0,\Gamma^{\mathcal A}}^h(0)]
 &\ge c_*M^{-1/2}.
 \label{eq:lower-antithetic}
\end{align}
Here $c_*$ is the constant defined in \eqref{eq:lower-constant}.
\end{theorem}
\begin{proof}
We first confirm the continuous-time value function is as claimed. We set
$\Phi(t,x)=x_1+1-\eu^{-(1-t)/2}\cos x_2$.
Differentiation gives
\[
 \partial_t\Phi+
 \tfrac12\tr\bigl(\sigma(u)\sigma(u)^\top D_x^2\Phi\bigr)=0
 \quad\text{for every }u\in\U,
 \qquad \Phi(1,\cdot)=g.
\]
Thus $\Phi$ is a classical solution of the Hamilton--Jacobi--Bellman equation with terminal condition $g$ and at most linear growth.
Since every action attains the Hamiltonian minimum,
any constant control satisfies the optimality condition
in \citep[Theorem~IV.3.1]{Fleming2006}, which gives $V=\Phi$.

For any sequence $\Lambda=(\lambda_1,\ldots,\lambda_N)$ of
noise laws on $\R^2$ with finite first moments, let
$A_k^\Lambda$ and $B_k^\Lambda$ be the scalar value functions
with terminal costs $x_1$ and $1-\cos x_2$ and Euler steps
\[
 \begin{aligned}
 F^{h,1}(x_1,u_1,z)
 &=x_1+(h/2)^{1/2}\,u_1(z_1+z_2),\\
 F^{h,2}(x_2,u_2,z)
 &=x_2+h^{1/2}\,[u_2z_1+(1-u_2)z_2].
 \end{aligned}
\]
Linearity of integration and the product action set give,
for $\phi(x_1,x_2)=a(x_1)+b(x_2)$,
\[
 [T_k^\Lambda\phi](x_1,x_2)
 =\min_{u_1\in\{-1,1\}}\int_{\R^2}
   a(F^{h,1}(x_1,u_1,z))\lambda_{k+1}(\mathrm d z)
 +\min_{u_2\in\{0,1\}}\int_{\R^2}
   b(F^{h,2}(x_2,u_2,z))\lambda_{k+1}(\mathrm d z).
\]
Backward induction from the additive terminal cost consequently implies
\begin{equation}\label{eq:lower-separation}
 v_k^\Lambda(x_1,x_2)
 =A_k^\Lambda(x_1)+B_k^\Lambda(x_2).
\end{equation}
\cref{ex:mean} and \eqref{eq:scalar-value} give
\[
 A_k^P(x_1)=x_1,\qquad
 B_k^P(x_2)=1-\rho_h^{N-k}\cos x_2.
\]
The sum \eqref{eq:lower-separation} with $\Lambda = P$ equals $\Phi(t_k,x)$, which completes the proof of \eqref{eq:lower-values}.

For ordinary sampling, we have from \cref{ex:mean} that
\(
 \E[A_0^P(0)-A_0^{\Gamma^{\mathcal O}}(0)]
 =\sqrt{2/\pi}\,(hM)^{-1/2}.
\)
Regarding
\(
 \E[B_0^P(0)-B_0^{\Gamma^{\mathcal O}}(0)],
\)
evaluating the constant policy
$u_2=0$ at every stage gives
\[
\mathbb E B_0^{\Gamma^{\mathcal O}}(0)
\le \mathbb E \int
\left(1-\cos\left(h^{1/2}\sum_{j=1}^N z_{j,2}\right)\right)
\prod_{j=1}^N \gamma_j^{\mathcal O}(dz_j)
=1-\rho_h^N=B_0^P(0),
\]
where the equality follows from independence of the stage
batches and unbiasedness of $\Gamma^{\mathcal O}$.
Combining these results through
\eqref{eq:lower-separation} proves
\eqref{eq:lower-ordinary}.

For antithetic sampling, the empirical mean of $(z_1+z_2)/\sqrt2$ is exactly zero,
so \cref{ex:mean} gives $A_k^{\Gamma^{\mathcal A}}(x_1)=x_1$.
Thus only  $\E [B_0^P(0) - B_0^{\Gamma^{\mathcal A}}(0)]$
contributes to
the left-hand side of \eqref{eq:lower-antithetic} with the right-hand side following from \eqref{eq:scalar-lower}.
\end{proof}

\paragraph{Interpretation of the sampling lower bounds.}
\Cref{thm:lower} establishes necessity of the ordinary and
antithetic sampling terms, not sharpness of the $h^{1/4}$ regularization
remainder or a lower bound for every centered rule.
In the cosine example, a rule with identical coordinate marginals, such as tensor-product Gauss--Hermite quadrature \citep[Section~3.2, equation~(3.9)]{Picarelli2020}, assigns the same approximate expectation to both actions and removes our source of selection bias. The approximation may still differ, however, from the Gaussian expectation.
We also note that the ordinary-sampling divergence under time refinement relies on the unbounded terminal cost: with bounded running and terminal costs,
\[
 \EE{\Gamma}{P}\le 2\bigl(\|g\|_\infty+T\|\ell\|_\infty\bigr).
\]
The error cannot then diverge, but it need not vanish, as the cosine example shows.

\subsection{Sharpness of the $M^{-1/m}$ rate}
\label{sec:global-sharpness}\label{app:unweighted-obstruction}

The next result shows that the sample-size exponent $1/m$ in
\cref{thm:transport} is sharp for $m\ge3$.
It reflects a geometric limitation of approximating the
Gaussian law by finitely supported measures.

\begin{proposition}[Sharpness of the $M^{-1/m}$ rate]
\label[proposition]{prop:unweighted-obstruction}
For each integer $m\ge3$, there is a problem satisfying
\textup{\cref{ass:data}}, with $T=h=1$, $n=m$, no control,
$f=\ell=0$, $\sigma=I_m$, and a bounded globally Lipschitz
terminal cost $g$, such that every probability integration
law $\gamma$ supported on at most $M$ points satisfies
\begin{equation}\label{eq:unweighted-obstruction}
 \|v_0^\gamma-v_0^P\|_\infty
 \ge c_mM^{-1/m},\qquad M\ge1.
\end{equation}
The problem and $c_m>0$ are independent of $M$ and $\gamma$.
In particular, the bound holds pathwise for all three Monte Carlo
constructions.
\end{proposition}

\begin{proof}
We first prove the lower bound over a fixed class of Lipschitz test
functions.

Let $\mathcal K=[-1,1]^m$, and let $\mathcal F$ be the class of
$1$-Lipschitz functions $\varphi:\R^m\to[0,1]$ supported on
$\mathcal K$. We fix a probability law $\gamma$ supported on
a set $A$ of at most $M$ points. The function
\[
 \varphi_A(z)
 =\min\{\operatorname{dist}(z,A),
         \operatorname{dist}(z,\mathcal K^c),1\}
\]
belongs to $\mathcal F$ and satisfies
$\int\varphi_A\,\dd\gamma=0$.
We set $\mathcal K_0=[-1/2,1/2]^m$, whose volume is one.
We choose $\kappa_m\in(0,1/2]$ such that
$\omega_m\kappa_m^m\le1/2$, where $\omega_m$ is the volume
of the Euclidean unit ball, and set $r=\kappa_mM^{-1/m}$.
The radius-$r$ balls centered at points of $A$ cover at most
half of $\mathcal K_0$. On the remaining part,
$\operatorname{dist}(z,A)\ge r$ and
$\operatorname{dist}(z,\mathcal K^c)\ge1/2$, so $\varphi_A(z)\ge r$.
The Gaussian density has a positive lower bound $p_m$ on $\mathcal K_0$.
Consequently,
\begin{equation}\label{eq:finite-support-test-bound}
 \sup_{\varphi\in\mathcal F}
 \Bigl|\int\varphi\,\dd(\gamma-P)\Bigr|
 \ge \int\varphi_A\,\dd P
 \ge c_mM^{-1/m},
 \qquad c_m=\tfrac12p_m\kappa_m.
\end{equation}

The test function $\varphi_A$ depends on the noise law.
To obtain a single control problem, we construct a fixed terminal
cost whose translates approximate every function in $\mathcal F$.

By the Arzelà--Ascoli theorem on $\mathcal K$, the class $\mathcal F$
has a countable uniformly dense subset $(\varphi_j)_{j\ge1}$.
We choose a sequence $(j_r)_{r\ge1}$ in which every positive integer
occurs infinitely often. We set $R_r=r+1$ and choose centers
$a_r$ successively along the first coordinate axis so that
\[
 |a_r-a_s|>R_r+R_s+2\sqrt m,
 \qquad r\ne s.
\]
We define
\[
 g(x)=\sup_{r\ge1}\varphi_{j_r}(x-a_r).
\]
As a supremum of functions taking values in $[0,1]$ with
Lipschitz constant one, $g$ is bounded and globally
$1$-Lipschitz. Since $\mathcal K\subset\{|z|\le\sqrt m\}$,
the separation of the centers ensures that
\[
 g(a_r+z)=\varphi_{j_r}(z)
 \qquad\text{whenever }|z|\le R_r.
\]
For any probability law $\gamma$ and each fixed $j$, we let
$r\to\infty$ along indices with $j_r=j$. Since both functions
take values in $[0,1]$,
\[
 \Bigl|
 \int\bigl[g(a_r+z)-\varphi_j(z)\bigr]
       \,(\gamma-P)(\mathrm d z)
 \Bigr|
 \le \gamma\{|z|>R_r\}+P\{|z|>R_r\}
 \rightarrow0.
\]
The supremum over initial states therefore dominates the integration
error of every $\varphi_j$. Using the uniform density of
$(\varphi_j)_{j\ge1}$ in $\mathcal F$ gives
\begin{equation}\label{eq:encoded-test-class}
 \sup_{x\in\R^m}
 \Bigl|\int g(x+z)\,(\gamma-P)(\mathrm d z)\Bigr|
 \ge \sup_{j\ge1}
       \Bigl|\int\varphi_j\,\dd(\gamma-P)\Bigr|
 =\sup_{\varphi\in\mathcal F}
       \Bigl|\int\varphi\,\dd(\gamma-P)\Bigr|.
\end{equation}

Finally, we take $T=h=1$, $n=m$, no control, $f=\ell=0$,
$\sigma=I_m$, and the terminal cost $g$ constructed above.
These data satisfy \cref{ass:data}, and the one-step values are
\[
 v_0^\gamma(x)=\int g(x+z)\,\gamma(\mathrm d z),
 \qquad
 v_0^P(x)=\int g(x+z)\,P(\mathrm d z).
\]
Combining \eqref{eq:encoded-test-class} with
\eqref{eq:finite-support-test-bound} proves
\eqref{eq:unweighted-obstruction}.
The terminal cost $g$ is fixed independently of $M$ and $\gamma$;
only the test function $\varphi_A$ depends on the
noise law.
\end{proof}

\section{Conclusion}\label{sec:conclusion}

The distinction between zero mean in expectation and zero
mean in every batch has consequences for time refinement.
If the population Euler error is $O(h^{1/4})$, our bounds
give sufficient choices $h\asymp\delta^4$ and
$M\asymp\delta^{-6}$ for ordinary Gaussian sampling,
compared with $M\asymp\delta^{-2}$ for antithetic or
recentered sampling, to obtain total error of order
$\delta$.
The consistency results show that
our $h^{1/4}$ regularization remainder is not an intrinsic fixed-mesh
sampling error.
Extending the analysis to approximate Bellman
updates and policy-performance guarantees is left for
future work.

\appendix
\section{Supporting proofs and calculations}\label{app:technical}

\subsection{Well-posedness and policies}\label{app:policies}

\begin{lemma}[Policy representation]
\label[lemma]{lem:wellposed}
Under \textup{\cref{ass:data}}, let
$\Lambda=(\lambda_1,\ldots,\lambda_N)$ be a sequence
of probability measures on $\R^m$ with finite first
moments. The functions $v_k^\Lambda$, $k=0,\ldots,N$,
are defined by the Bellman recursion \eqref{eq:recursion}.
\begin{enumerate}[label=\textup{(\alph*)},ref=(\alph*),leftmargin=*,nosep]
\item\label{item:wellposed-regularity}
These functions are finite, globally Lipschitz,
and of at most linear growth.
The constants may depend on $h$ and $\Lambda$.

\item\label{item:wellposed-policies}
For a starting stage $k$ and state $x$, let
$(\mathcal H_j)_{j=k}^N$ be a filtration.
For $j=k,\ldots,N-1$, suppose that $Z_{j+1}$ is
$\mathcal H_{j+1}$-measurable and satisfies
\eqref{eq:conditional-law}.
Then $v_k^\Lambda(x)$ equals the infimum in
\eqref{eq:discrete-cost}, taken over controls
$u=(u_j)_{j=k}^{N-1}$ for which each action $u_j$
is $\mathcal H_j$-measurable and takes values in $\U$.
The infimum is unchanged when restricted to Markov
policies $u_j=\pi_j(S_j)$ with Borel maps
$\pi_j:\R^n\to\U$.
\end{enumerate}
\end{lemma}

\begin{proof}[Proof sketch]
For \ref{item:wellposed-regularity}, let $\Phi_j(x,u)$ denote the
objective in \eqref{eq:bellman} with $\phi=v_{j+1}^\Lambda$.
The growth conditions in \cref{ass:data} give
\begin{equation}\label{eq:transition-growth}
 |F_j^h(x,u,z)|
 \le (1+hL_1)|x|+hL_1+h^{1/2}L_1|z|.
\end{equation}
If $v_{j+1}^\Lambda$ is globally Lipschitz with at most linear
growth, the data bounds and the finite first moment of $\lambda_{j+1}$
give both a linear growth bound and a state Lipschitz bound for
$\Phi_j(\cdot,u)$, uniformly in $u$. Starting from $v_N^\Lambda=g$,
backward induction using \eqref{eq:infimum-comparison} gives
\ref{item:wellposed-regularity}.

For \ref{item:wellposed-policies}, \eqref{eq:transition-growth}
gives $\E|S_j|<\infty$ by forward induction, so the costs and value
functions along the trajectory are integrable. By \eqref{eq:conditional-law},
$\Phi_j(S_j,u_j)$ equals $h\ell(t_j,S_j,u_j)$ plus
$\E[v_{j+1}^\Lambda(S_{j+1})\mid\mathcal H_j]$.
Taking expectations in $v_j^\Lambda(S_j)\le\Phi_j(S_j,u_j)$
and summing over $j=k,\ldots,N-1$ shows that $v_k^\Lambda(x)$
is a lower bound on every admissible cost.

For the upper bound, continuity of $\Phi_j(x,\cdot)$ follows
by dominated convergence using \eqref{eq:transition-growth}.
Given $\eta>0$, let $\pi_j(x)$ be the first action $u$ in a fixed
countable dense sequence in $\U$ with
$\Phi_j(x,u)\le v_j^\Lambda(x)+\eta/N$.
Such an action exists by continuity, and $\pi_j$ is Borel because each
comparison is between Borel functions of $x$. Under this Markov policy,
taking expectations and summing gives a cost of at most
$v_k^\Lambda(x)+\eta$. Letting $\eta\downarrow0$ proves
\ref{item:wellposed-policies}.
\end{proof}

\subsection{Gaussian sampling bounds}\label{app:moments}

\begin{proof}[Proof of \textup{\cref{cor:samplers}}]
Let $Z_1,\ldots,Z_M$ be independent copies of
$Z\sim P$. We write $\bar Z=M^{-1}\sum_{i=1}^{M} Z_i$.

We first verify \cref{ass:sampling} for ordinary and
recentered sampling. Both have the same centered points
$Z_i-\bar Z$, and hence the same moments $Q,H,R$.
From the definition of $Q$,
\[
 Q-I_m=\frac1M\sum_{i=1}^{M}(Z_iZ_i^\top-I_m)
       -\bar Z\bar Z^\top.
\]
The summands in the first term are independent and
have mean zero. Since
$\|\bar Z\bar Z^\top\|_{\F}^2=|\bar Z|^4$
and
$\bar Z$ is equal in distribution to $M^{-1/2}Z$,
we obtain
\[
 \E\|Q-I_m\|_{\F}^2
 \le \frac2M\E\|ZZ^\top-I_m\|_{\F}^2
      +2\E|\bar Z|^4
 =\frac2M\E\|ZZ^\top-I_m\|_{\F}^2
      +\frac2{M^2}\E|Z|^4
 \le \frac{C_m}{M}.
\]
Expanding the centered third moment $H$ and applying
the triangle inequality yields
\[
 \|H\|_{\F}
 \le \biggl\|\frac1M\sum_{i=1}^{M}Z_i^{\otimes3}\biggr\|_{\F}
 +3|\bar Z|\frac1M\sum_{i=1}^{M}|Z_i|^2
 +2|\bar Z|^3.
\]
Since $\E Z^{\otimes3}=0$, we can use Jensen's inequality
and independence to bound the first term:
\[
\E\biggl\|\frac1M\sum_{i=1}^{M} Z_i^{\otimes3}\biggr\|_{\F}
\le
\Biggl(\E\biggl\|\frac1M\sum_{i=1}^{M} Z_i^{\otimes3}
\biggr\|_{\F}^2\Biggr)^{1/2}
=\biggl(\frac1{M^2}\sum_{i=1}^{M}
\E\|Z_i^{\otimes3}\|_{\F}^2\biggr)^{1/2}
=\left(\frac{\E|Z|^6}{M}\right)^{1/2}.
\]
For the middle term, the Cauchy--Schwarz inequality
shows that
\[
 \E\biggl[|\bar Z|\frac1M\sum_{i=1}^{M}|Z_i|^2\biggr]
 \le (\E|\bar Z|^2)^{1/2}
 \Biggl[\E\biggl(\frac1M\sum_{i=1}^{M}|Z_i|^2\biggr)^2
 \Biggr]^{1/2}
 \le \left(\frac{m}{M}\E|Z|^4\right)^{1/2},
\]
where the last inequality follows from Jensen's
inequality and $\E|\bar Z|^2=m/M$.
Finally, $\E|\bar Z|^3=M^{-3/2}\E|Z|^3$.
Combining these estimates gives
\[
 \E\|H\|_{\F}\le C_mM^{-1/2}.
\]
Also, each $Z_i-\bar Z$ has law
$\mathcal N(0,(1-M^{-1})I_m)$, so
\[
 \E R=(1-M^{-1})^2\E|Z|^4\le C_m.
\]

For antithetic sampling, $M=2\nu$ and $H=0$.
The covariance $Q$ is the average of $Z_iZ_i^\top$
over $\nu$ independent Gaussian generators, and
$R$ is the corresponding average of $|Z_i|^4$.
Thus
\[
 \E\|Q-I_m\|_{\F}^2
 =\frac1\nu\E\|ZZ^\top-I_m\|_{\F}^2
 \le\frac{C_m}{M},
 \qquad
 \E R=\E|Z|^4\le C_m.
\]
All three schemes therefore satisfy
\cref{ass:sampling} with $K$ depending only on $m$.

For ordinary sampling, $b_j$ has the law of
$\bar Z$, and hence
\[
 h^{1/2}\sum_{j=1}^N\E|b_j|
 =T\E|Z|(hM)^{-1/2}.
\]
For antithetic and recentered sampling, $b_j=0$.
For ordinary sampling, the term $M^{-1/2}$ in
\cref{thm:main} is bounded by $(hM)^{-1/2}$ because
$h\le1$. Substituting the preceding estimates into
that theorem completes the proof.
\end{proof}

\subsection{A sufficient population approximation result}
\label{app:population}

\begin{proposition}[Population Euler approximation]
\label[proposition]{prop:population}
Under \textup{\cref{ass:data}}, suppose that $\U$ is compact,
$f$ and $\ell$ are uniformly bounded, and
$f,\sigma,\ell$ are uniformly $1/2$-H\"older in time,
independently of $(x,u)$.
Then the piecewise-constant extension $V^h$ of the
population Euler values $V_k^h$, defined in
\textup{\cref{sec:continuous}}, satisfies
\begin{equation*}
 \sup_{t\in[0,T],\,x\in\R^n}
 |V^h(t,x)-V(t,x)|\le Ch^{1/4}.
\end{equation*}
Here $C$ may also depend on the additional uniform
bounds and H\"older constants in the proposition.
\end{proposition}

\begin{proof}
The Gaussian increment $h^{1/2}Z$, $Z\sim P$,
satisfies the moment conditions in
\citep[Section~3.2]{Jakobsen2019}.
After reversing signs for minimization, the approximation
estimate there implies
\[
 \max_{0\le k\le N}\sup_{x\in\R^n}
 |V_k^h(x)-V(t_k,x)|\le Ch^{1/4}.
\]

To extend this estimate between mesh points, we use the
time-regularity argument in \citep[proof of Corollary~3.4(i)]{Krylov1999},
which gives
\[
 |V(t,x)-V(s,x)|\le C|t-s|^{1/2},
 \qquad s,t\in[0,T],\quad x\in\R^n.
\]
In our undiscounted setting, the argument uses only
the Lipschitz continuity of $g$, so its boundedness
is unnecessary.

For $t\in[t_k,t_{k+1})$, we have
$V^h(t,x)=V_k^h(x)$ and $0\le t-t_k<h$. Hence
\[
 |V^h(t,x)-V(t,x)|
 \le |V_k^h(x)-V(t_k,x)|
       +|V(t_k,x)-V(t,x)|
 \le C(h^{1/4}+h^{1/2})
 \le Ch^{1/4},
\]
where the last inequality uses $h\le1$.
Since $V^h(T,\cdot)=V(T,\cdot)=g$, the proof is complete.
\end{proof}

\subsection{Normalization of a centered law}
\label{app:normalization}

\begin{proof}[Proof of \textup{\cref{lem:normalization}}]
We set $d=\|Q-I_m\|_{\F}$. Since $\|Q\|_{\op}\le1+d$,
we have $0\le a^2-1\le d$ from $a\ge1$.
Scaling by $a^{-1}$ gives covariance $\widetilde Q$,
$0\preceq\widetilde Q\preceq I_m$ and $\tr Q\le ma^2$.
Coupling $Z\sim\mu$ with $Z/a$ yields
\[
 W_2(\mu,\widetilde\mu)
 \le(1-a^{-1})\sqrt{\tr Q}
 \le\sqrt m\,(a-1)
 \le\sqrt m\,d.
\]
For $X\sim P$, we couple $X$ with $\widetilde Q^{1/2}X$.
Applying $1-\sqrt t\le1-t$ to the eigenvalues of
$\widetilde Q$ and then applying the triangle inequality
to $(\widetilde Q - a^{-2}I_m) + (a^{-2} - 1)I_m$ gives
\[
 W_2(G_Q,P)
 \le\|\widetilde Q^{1/2}-I_m\|_{\F}
 \le\|\widetilde Q-I_m\|_{\F}
 \le a^{-2}\bigl[d+\sqrt m\,(a^2-1)\bigr]
 \le(1+\sqrt m)d.
\]
Adding the displays proves \eqref{eq:normalization-w2}.
The moment relations for $\widetilde\mu$ follow by scaling.
The Gaussian fourth-moment formula gives
\[
 m_4(G_Q)
 =(\tr\widetilde Q)^2+2\tr(\widetilde Q^2)
 \le m(m+2).
\]
All calculations remain valid when $Q$ is singular.
\end{proof}

\subsection{The cosine moment estimates}\label{app:cosine}
Here we verify \eqref{eq:D-moments}. Let $Z_1,Z_2$ be independent standard normal random variables and
$Y=\cos(h^{1/2}Z_1)-\cos(h^{1/2}Z_2)$. From \eqref{eq:cosine-gaussian} and the identity $\cos^2a=(1+\cos2a)/2$,
\[
 \E Y=0,\qquad
 \E Y^2=2\left[\tfrac12(1+\eu^{-2h})-\eu^{-h}\right]=(1-\eu^{-h})^2.
\]
Also $|1-\cos a|\le a^2/2$, which implies
\[
 |Y|\le\frac{h}{2}(Z_1^2+Z_2^2),\qquad
 \E Y^4\le\frac{h^4}{16}\E(Z_1^2+Z_2^2)^4=24h^4.
\]
The last expectation equals $384$, obtained by expanding the fourth power and using the Gaussian even moments.
In \cref{lem:cosine}, $\Xi=\nu^{-1}\sum_{i=1}^{\nu}Y_i$ where $Y_1,\ldots,Y_\nu$ are independent copies of $Y$.
In expanding the fourth power of $\Xi$, only four equal indices or two pairs of equal indices survive. Hence
\[
 \E \Xi^2=\frac{\E Y^2}{\nu},\qquad
 \E \Xi^4=\frac{\nu\E Y^4+3\nu(\nu-1)(\E Y^2)^2}{\nu^4}
 \le\frac{h^4}{\nu^2}\left(\frac{24}{\nu}+3\right)
 \le\frac{27h^4}{\nu^2}.
\]

\section*{Acknowledgments}
GPT 5.6 and 6 assisted with literature exploration, mathematical development and checking, and exposition during preparation of this manuscript. The authors remain responsible for the manuscript's contents.\footnote{%
This article has been authored by an employee of National Technology \&
Engineering Solutions of Sandia, LLC under Contract No. DE-NA0003525 with the
U.S. Department of Energy (DOE).
The employee owns all right, title and interest in and to the article and is
solely responsible for its contents.
The United States Government retains and the publisher, by accepting the
article for publication, acknowledges that the United States Government
retains a non-exclusive, paid-up, irrevocable, world-wide license to publish or
reproduce the published form of this article or allow others to do so, for
United States Government purposes.
The DOE will provide public access to these results of federally sponsored
research in accordance with the DOE Public Access Plan:
\url{https://www.energy.gov/downloads/doe-public-access-plan}.
This paper describes objective technical results and analysis.
Any subjective views or opinions that might be expressed in the paper do not
necessarily represent the views of the U.S. Department of Energy or the United
States Government.}

\bibliography{sampling-error-time-discretized-optimal-control.bbl}

\begin{thebibliography}{10}

\bibitem{Bonnans2019}
J.~F. Bonnans, J. Gianatti, and F.~J. Silva.
\newblock On the time discretization of stochastic optimal control problems:
  {T}he dynamic programming approach.
\newblock {\em ESAIM: Control, Optimisation and Calculus of Variations}, 25:63,
  2019.
\newblock \href {https://doi.org/10.1051/cocv/2018045}
  {\path{doi:10.1051/cocv/2018045}}.

\bibitem{Chen2026}
X. Chen, E.~S. Isik, and J. Milz.
\newblock Asymptotic analysis of empirical dynamic programming in
  infinite-horizon stochastic optimal control.
\newblock Preprint, 2026.
\newblock \href {http://arxiv.org/abs/2607.21520v1}
  {\path{arXiv:2607.21520v1}}.

\bibitem{Fahim2011}
A. Fahim, N. Touzi, and X. Warin.
\newblock A probabilistic numerical method for fully nonlinear parabolic
  {PDEs}.
\newblock {\em Ann. Appl. Probab.}, 21(4):1322--1364, 2011.
\newblock \href {https://doi.org/10.1214/10-AAP723}
  {\path{doi:10.1214/10-AAP723}}.

\bibitem{Fleming2006}
W.~H. Fleming and H.~M. Soner.
\newblock {\em Controlled {Markov} processes and viscosity solutions},
  volume~25 of {\em Appl. Math.}
\newblock New York: Springer, 2006.
\newblock \href {https://doi.org/10.1007/0-387-31071-1}
  {\path{doi:10.1007/0-387-31071-1}}.

\bibitem{Fournier2015}
N. Fournier and A. Guillin.
\newblock On the rate of convergence in {Wasserstein} distance of the empirical
  measure.
\newblock {\em Probab. Theory Related Fields}, 162(3--4):707--738, 2015.
\newblock \href {https://doi.org/10.1007/s00440-014-0583-7}
  {\path{doi:10.1007/s00440-014-0583-7}}.

\bibitem{Fuellner2025}
C. F\"ullner and S. Rebennack.
\newblock Stochastic dual dynamic programming and its variants: {A} review.
\newblock {\em SIAM Rev.}, 67(3):415--539, 2025.
\newblock \href {https://doi.org/10.1137/23M1575093}
  {\path{doi:10.1137/23M1575093}}.

\bibitem{Haskell2016}
W.~B. Haskell, R. Jain, and D. Kalathil.
\newblock Empirical dynamic programming.
\newblock {\em Math. Oper. Res.}, 41(2):402--429, 2016.
\newblock \href {https://doi.org/10.1287/moor.2015.0733}
  {\path{doi:10.1287/moor.2015.0733}}.

\bibitem{Jakobsen2019}
E.~R. Jakobsen, A. Picarelli, and C. Reisinger.
\newblock Improved order {$1/4$} convergence for piecewise constant policy
  approximation of stochastic control problems.
\newblock {\em Electron. Commun. Probab.}, 24:Paper No. 59, 10, 2019.
\newblock \href {https://doi.org/10.1214/19-ecp256}
  {\path{doi:10.1214/19-ecp256}}.

\bibitem{Kleywegt2002}
A.~J. Kleywegt, A. Shapiro, and T. Homem-de Mello.
\newblock The sample average approximation method for stochastic discrete
  optimization.
\newblock {\em SIAM J. Optim.}, 12(2):479--502, 2002.
\newblock \href {https://doi.org/10.1137/S1052623499363220}
  {\path{doi:10.1137/S1052623499363220}}.

\bibitem{Krylov1999}
N.~V. Krylov.
\newblock Approximating value functions for controlled degenerate diffusion
  processes by using piece-wise constant policies.
\newblock {\em Electron. J. Probab.}, 4:19, 1999.
\newblock Id/No 2.
\newblock \href {https://doi.org/10.1214/EJP.v4-39}
  {\path{doi:10.1214/EJP.v4-39}}.

\bibitem{Krylov2000}
N.~V. Krylov.
\newblock On the rate of convergence of finite-difference approximations for
  {Bellman's} equations with variable coefficients.
\newblock {\em Probab. Theory Related Fields}, 117(1):1--16, 2000.
\newblock \href {https://doi.org/10.1007/s004400050264}
  {\path{doi:10.1007/s004400050264}}.

\bibitem{Kushner2001}
H.~J. Kushner and P. Dupuis.
\newblock {\em Numerical methods for stochastic control problems in continuous
  time}, volume~24 of {\em Applications of Mathematics}.
\newblock Springer, New York, second edition, 2001.
\newblock Stochastic Modelling and Applied Probability.
\newblock \href {https://doi.org/10.1007/978-1-4613-0007-6}
  {\path{doi:10.1007/978-1-4613-0007-6}}.

\bibitem{Lan2024}
G. Lan and A. Shapiro.
\newblock Numerical methods for convex multistage stochastic optimization.
\newblock {\em Foundations and Trends in Optimization}, 6(2):63--144, 05 2024.
\newblock \href {https://doi.org/10.1561/2400000044}
  {\path{doi:10.1561/2400000044}}.

\bibitem{LedouxZhu2021}
M. Ledoux and J.-X. Zhu.
\newblock On optimal matching of {Gaussian} samples {III}.
\newblock {\em Probab. Math. Statist.}, 41(2):237--265, 2021.
\newblock \href {https://doi.org/10.37190/0208-4147.41.2.3}
  {\path{doi:10.37190/0208-4147.41.2.3}}.

\bibitem{Milz2025CLT}
J. Milz and A. Shapiro.
\newblock Central limit theorems for sample average approximations in
  stochastic optimal control.
\newblock Preprint, 2025.
\newblock \href {http://arxiv.org/abs/2508.01942v2}
  {\path{arXiv:2508.01942v2}}, \href
  {https://doi.org/10.48550/arxiv.2508.01942}
  {\path{doi:10.48550/arxiv.2508.01942}}.

\bibitem{Pagnoncelli2024}
B.~K. Pagnoncelli, T. Homem-de Mello, G. Lagos, P. Casta{\~n}eda, and J.
  Garc{\'i}a.
\newblock Solving constrained consumption--investment problems by decomposition
  algorithms.
\newblock {\em European J. Oper. Res.}, 319(1):292--302, 2024.
\newblock \href {https://doi.org/10.1016/j.ejor.2024.06.027}
  {\path{doi:10.1016/j.ejor.2024.06.027}}.

\bibitem{Picarelli2020}
A. Picarelli and C. Reisinger.
\newblock Probabilistic error analysis for some approximation schemes to
  optimal control problems.
\newblock {\em Systems Control Lett.}, 137:104619, 2020.
\newblock \href {https://doi.org/10.1016/j.sysconle.2019.104619}
  {\path{doi:10.1016/j.sysconle.2019.104619}}.

\bibitem{Tran2022}
K.~Q. Tran, B.~T.~N. Le, and G. Yin.
\newblock Harvesting of a stochastic population under a mixed regular-singular
  control formulation.
\newblock {\em J. Optim. Theory Appl.}, 195:1106--1132, 2022.
\newblock \href {https://doi.org/10.1007/s10957-022-02127-7}
  {\path{doi:10.1007/s10957-022-02127-7}}.

\bibitem{Villani2009}
C. Villani.
\newblock {\em Optimal Transport: Old and New}, volume 338 of {\em Grundlehren
  der mathematischen Wissenschaften}.
\newblock Springer, Berlin, Heidelberg, 2009.
\newblock \href {https://doi.org/10.1007/978-3-540-71050-9}
  {\path{doi:10.1007/978-3-540-71050-9}}.

\bibitem{Zhu2022}
D. Zhu and H. Zheng.
\newblock Effective approximation methods for constrained utility maximization
  with drift uncertainty.
\newblock {\em J. Optim. Theory Appl.}, 194:191--219, 2022.
\newblock \href {https://doi.org/10.1007/s10957-022-02015-0}
  {\path{doi:10.1007/s10957-022-02015-0}}.

\end{thebibliography}
\end{document}